\documentstyle[amscd,amssymb,verbatim,12pt]{amsart}
\input xypic \xyoption{curve} 

\def\hcorrection#1{\advance\hoffset by #1 }
\def\vcorrection#1{\advance\voffset by #1 }

\vcorrection{-.7in}
\hcorrection{-.5in}

\newcommand{\B}[1]{{\bold#1}} 
\newcommand{\C}[1]{{\cal#1}} 

\newcommand{\one}{{\bf 1}}

\newcommand{\To}{\longrightarrow}
\newcommand{\sst}[1]{\scriptstyle{#1}}

\newcommand{\Gr}{$\tilde{\C{G}}r$}
\newcommand{\tens}[1]{{\otimes_\C{#1}}}

\theoremstyle{plain}
\newtheorem{th}{Theorem}[section]
\newtheorem{cor}{Corollary}[section]

\newtheorem{lem}{Lemma}[section]
\newtheorem{prop}{Proposition}[section]

\theoremstyle{definition}
\newtheorem{defin}{Definition}[section]

\theoremstyle{definition}
\newtheorem{example}{Example}[section]

\newtheorem{rem}{Remark}[section]
\newtheorem{notation}{Notation}[section]

\numberwithin{equation}{section}

\begin{document}
\pagestyle{plain}
\addtolength{\footskip}{.3in}

\title{On Parity Complexes\\ and \\Non-abelian cohomology}
\author{Lucian M. Ionescu}
\address{Mathematics Department\\Kansas State University\\
             Manhattan, Kansas 66502}
\email{luciani@@math.ksu.edu}
\date{August 10, 1998}

\begin{abstract}
To characterize categorical constraints - associativity, commutativity
and monoidality - in the context of quasimonoidal categories,
from a cohomological point of view,
we define the notion of a parity (quasi)complex.

Applied to groups
gives non-abelian cohomology. The categorification - functor
from groups to monoidal categories - provides the correspondence
between the respective parity (quasi)complexes and allows to interpret
1-cochains as functors, 2-cocycles - monoidal structures,
3-cocycles - associators.

The cohomology spaces $H^3, H^2, H^1, H^0$ correspond as usual to
quasi-extensions, extensions, split extensions and invariants,
as in the abelian case.

A larger class of commutativity constraints for monoidal categories
is identified. It is naturally associated with 
coboundary Hopf algebras.
\end{abstract}

\maketitle
\tableofcontents


\section{Introduction}\label{S:intro}
The cohomological deformation theory for braided monoidal categories
was provided in \cite{CY1,CY2}. As remarked in \cite{Yet1}, section 3,
the coherence conditions are ``formally a cohomological condition
written multiplicatively'' as was shown in \cite{Dav1}. In both approaches
a fixed associator (``base point") is needed to ensure that the
composition of the functorial morphisms considered in a deformation of the
associativity constraint or of the monoidal structure of a monoidal functor,
is well defined.

Formal deformations of associators in a monoidal category are controlled
by a complex analogous to the Hochschild complex. It is a description
at the infinitesimal level.
To characterize from a cohomological point of view and in a global way
the basic algebraic constraints at a categorical level - associativity,
commutativity  and morphism condition - the concept of
a {\em parity complex} is introduced (definition \ref{Def:pqc}).
It provides a multiplicative analogue of the Hochschild cohomology. It is
based on the natural idea of ``sign separation'', which is appropriate
in the presence of non-commutativity or lack of inverses.

The standard parity complex associated to a group acting on a
non-abelian group provides the correct non-abelian cohomology in
low dimensions, as explained in section \ref{S:nacoh}.
Writing cohomological conditions in a direct manner - 
different versions in dimensions higher then two -
is at least not appealing \cite{Dedecker}, even in a categorical approach
\cite{Breen}.
In our approach non-comutativity (non-associativity) is thought of
alternatively as yielding curvature (holonomy) or monoidality.

Non-abelian cohomology characterizes group extensions as explained in
section \ref{SS:ext}. For arbitrary 2-cochains $(f,L)\in C^2(G,N)$
an associated quasi-extension (monoidal category) is associated
(theorem \ref{T:qext}).
The concept of {\em integrability} is introduced (definition \ref{Def:int}).
It corresponds both to the lack of ``curvature'' for the quasicomplex
($\delta^2=1$) and to the
existence of an extension of the group $G$ by a subgroup of the coefficient
group $N$. It is controlled by
an equation in dimension two (definition \ref{Def:int} and
equation \ref{E:MC}). It is shown that for
$0\le p\le 2$,  $p$@-cocycles are integrable (theorem \ref{T:ext}).


A group also defines a monoidal category, by a
functorial process we introduce in section \ref{S:categ},
and called categorification. It is a special instance of categorifying
a group extension (theorem \ref{T:categ}), thought of as an affine bundle.
It is different from the procedure defined in \cite{CY2}. It provides
the natural correspondence under which non-abelian cohomology
characterizes monoidal functors (proposition \ref{P:fiber}).
The integrability equation $(MC)$ is again a monoidality condition
(theorem \ref{T:MC}). The correspondence gives
the following ``dictionary'':
$$\begin{matrix}
& \quad  & \quad \underline{Groups} \quad &<->& \quad \underline{Categories}\quad\\
1-cochain &        & function                & & monoidal\ functor\\
1-cocycle &        & morphism                & & strict\ monoidal\ functor\\
2-cochain &        & quasiextension          & & monoidal category\\
2-cocycle &        & extension               & & strict \ monoidal\ category\\
3-cocycle &        & obstruction             & & (coherent)\ associator
\end{matrix}$$
Returning to the categorical level, to a quasimonoidal category there is
an associated parity complex in low dimensions
(diagram \ref{D:monpqc}). The approach slightly generalizes some of the
corresponding results from \cite{Dav1} (see also \cite{Sa}).
The coherence condition for a quasi-associator
is the 3-cocycle condition.

To study the compatibility
of a functor between two categories with the corresponding
products (not necessarily monoidal),
the alternate view point to a parity quasicomplex approach is in terms
of a biaction, as explained in section \ref{S:cohmc}. When specialized
to monoidal categories one obtains the characterization of the monoidal
structures of the functor. When the functor is the identity functor,
the characterization of monoidal structures of the category is
obtained, and finally that the {\em commutativity} constraint,
which naturally corresponds to categories of representations of
coboundary Hopf algebras (\cite{CP}), is proved as being
the monoidal structure condition for the identity functor between a
monoidal category and its opposite
(section \ref{SS:cc}; see also \cite{I}), i.e. the commutativity
constraints are the 2-cocycles. A braiding is a monoidal structure
of the ``identity'' functor from the underlying monoidal category
to its opposite (\cite{JS1,I}), being a special case.



\section{\bf Parity Quasicomplexes}\label{S:pqc}
Comparing the coherence condition for the associativity constraint
$\alpha$ of a monoidal category $(\C{C}, \otimes, \alpha)$:
%
%
\begin{gather}\diagram
 & \sst{(X\otimes Y)\otimes(Z\otimes W)}
  \drto^{\alpha_{X, Y, Z\otimes W}} & \\
\sst{((X \otimes Y) \otimes Z) \otimes W}
  \urto^{\alpha_{X\otimes Y,Z,W}} \dto_{\alpha_{X,Y,Z}\otimes I_W}
& & \sst{X\otimes(Y\otimes(Z\otimes W))}\\
\sst{(X\otimes(Y\otimes Z))\otimes W} \rrto^{\alpha_{X,Y\otimes Z,W}} &
& \sst{X\otimes((Y\otimes Z)\otimes W)} \uto_{I_X\otimes \alpha_{Y,Z,W}}
\enddiagram\label{D:pen}\end{gather}
with the Hochschild differential of the complex $(C^\bullet(A), d)$
associated to an associative algebra $(A,\mu)$, in the appropriate
degree:
\begin{equation}\label{E:pent}
d\alpha(x,y,z,w)=x \alpha(y,z,w)-\alpha(xy,z,w)+\alpha(x,yz,w)-
\alpha(x,y,zw)+\alpha(x,y,z)w
\end{equation}
it is clear that there is a close connection. The first is a
multiplicative analogue of the second.

The infinitesimal coherence condition, appearing in the formal deformation
of an associator, was studied by L. Crane and D. Yetter in \cite{CY1,Yet1}.

Diagram \ref{D:pen} is thought of as the non-linear version.

Recalling Gerstenhaber's comp operation
$f\circ_i g=f(\ ,...,g,...,\ ,\ )$ (\cite{GS}):
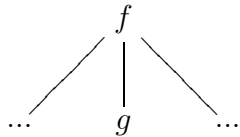
$$\diagram
& f \ar@{-}[dl] \ar@{-}[d] \ar@{-}[dr]\\
... & g & ...
\enddiagram$$
then $\alpha(xy,z,w)$ would be just $\alpha_1=\alpha\circ_1\otimes$, etc.
The diagram then reads:
\begin{gather}\diagram
 & \sst{(X\otimes Y)\otimes(Z\otimes W)}
  \drto^{\alpha_3} & \\
\sst{((X \otimes Y) \otimes Z) \otimes W}
  \urto^{\alpha_1} \dto_{\alpha_4}
& & \sst{X\otimes(Y\otimes(Z\otimes W))}\\
\sst{(X\otimes(Y\otimes Z))\otimes W} \rrto^{\alpha_2} &
& \sst{X\otimes((Y\otimes Z)\otimes W)} \uto_{\alpha_0}
\enddiagram\label{D:pent}\end{gather}

%
In passing from abelian to non-abelian cohomology a somewhat heuristic
approach is
to separate the pluses and minuses in the definition of the differential,
and consider pairs of (co)boundary operators $\partial^+, \partial^-$.
The approach avoids the completion procedure when dealing with ``rigs'',
fusion rings, etc. when opposites or inverses are not available.
The term is similar to that in \cite{St2}, where parity complexes are defined
as combinatorial structures.
The idea of odd faces being associated
with the source and even faces with the target is present.
The corresponding cohomological conditions in low dimesions
are also implicite in \cite{St1} (p.288).
\begin{defin}\label{Def:pqc}
A {\em parity quasicomplex} in a category $\C{C}$ is a sequence of
pairs of morphisms $(C_n,\partial_n^+,\partial_n^-)$
having equalizers for every $n$. 
If $\C{C}$ is additive, it is called a {\em parity complex} if
$(C_n,\partial_n=\partial_n^+-\partial_n^-)$ is a complex.

If $\C{C}$ is a concrete category, elements $f\in C^n$
such that $\partial^+f=\partial^-f$ are called
{\em cocycles} (or {\em first order  cocycles}).

If $C^n$ are monoids, $c$ and $c'\in C^n$ are {\em cobordant}
and  denoted as $c\overset{f}{\to}c'$
if there is an element $f\in C^{n-1}$,
called {\em cobordism} such that $\partial^-f=c$ and $\partial^+f=c'$
({\em closed cobordism} if $c=c'$).

If $c\circ\partial^-f=(\partial^+f)\circ c'$ then $c$ and $c'$ are
{\em cohomologous}, and denoted as $c\overset{f}{\sim}c'$
(or {\em second order cocycle} if $c=c'$).
\end{defin}
Recall that the equalizer (``the difference kernel'' \cite{ML})
of a pair of morphisms in an abelian category is just the kernel of
their difference and in category $\C{S}et$ is $\{x|f(x)=g(x)\}$.

The motivation for the
alternate terminology comes from the ``deformation equation'':
$$ (I+h \partial^+f)\circ(I+h x)=(I+h x)\circ(I+h \partial^-h)$$
It is satisfied to corresponding orders in $h$ if:
\begin{align}
\partial^+f+x&=\partial^-f+x, \qquad \text{first order}\notag\\
(\partial^+f)\circ x&=x\circ\partial^-f, \qquad \text{second order}
\end{align}
Some simple consequences are stated next:
\begin{prop}
If $(C_n,\partial_n^+,\partial_n^-)$ is a parity quasicomplex of
groups, then:\\
(i) A cocycle is a closed cobordism.\\
(ii) Two cobordant elements are cohomologous.\\
(iii) $\sim$ is an equivalence relation.\\
(iv) The cocycles form a group $Z^n$.\\
(v) If the groups $C^n$ are abelian, then being cobordant, cohomologous or
differing by a coboundary are equivalent statements. First and second
order cocycles coincide.\\
(vi) If $C^n$ are commutative rings, then a multiplicative non-zero
$2^{nd}$ order cocycle is a zero divisor of the corresponding coboundary.
\end{prop}
\begin{pf}
(i) Obviously any $f\in C^n$ is a cobordism
$\partial^-f\overset{f}{\to}\partial^+f$, and the closed ones are
precisely the cocycles.\\
(ii) If $c=\partial^-f$ and $c'=\partial^+f$, then clearly
$\partial^+f \circ c=c'\circ\partial^-f$.\\
(iii) For the group identity $e$, $\partial^\pm e=e$ and the relation
is reflexive. If $c\overset{f}{\sim}c'$ then $c\overset{f^{-1}}{\sim}c'$.
If the diagram commutes:
\xycompileto{art1}|{
\begin{align}\diagram
\dto_{a} \rto^{\partial^-c} & \dto^{a'} \rto^{\partial^-c} & \dto^{a''} \\
\rto^{\partial^+c} & \rto^{\partial^+c'} & 
\enddiagram\end{align}}
then $a\overset{c\circ c'}{\sim} a''$.\\
(v) Cobordant elements $c=\partial^+f$ and $c'=\partial^-f$ differ
by a coboundary $c-c'=\partial f$. Equivalently
$\partial^+f+c=c'+\partial^-f$.

A second order cocycle $\partial^+f+c=c+\partial^-f$ is just a
cocycle $\partial f=0$.\\
(vi) A multiplicative $2^{nd}$ order cocycle
$\partial^+f\circ c=c\circ\partial^-f$ is a zero divisor
$(\partial f)\circ c=0$.
\end{pf}
Comparing with a complex, the substitute for kernel is the equalizer,
and the analog for the coimage is the quotient by the $\sim$ relation.
The {\em cohomology spaces} $H^n$ are the quotient spaces $Z^n/\sim$.

As a special case, we have:
\begin{lem}
If a parity quasicomplex of groups is {\em central}
$Im\partial^-_{n-1}\subset Cen(C^n)$, then
$B^n=\{(\partial^+c)(\partial^-c)^{-1}|c\in C^{n-1}\}\cap Z^n$
is a group and $H^n\cong Z^n/B^n$
\end{lem}
\begin{pf}
Note that:
$$(\partial^+a)(\partial^-a)^{-1}(\partial^+b)(\partial^-b)^{-1}=
(\partial^+a)(\partial^+b)(\partial^-b)^{-1}(\partial^-a)^{-1}=
(\partial^+(ab))(\partial^-(ab))^{-1}$$
Also $c\overset{f}{\sim}c'$, or $\partial^+f c=c' \partial^-f$, is equivalent
to $c'=(\partial^+f) (\partial^-f)^{-1}) c$. Thus $\sim$ equivalence
classes are the left cosets of $B^n$.
\end{pf}
If $(C^\bullet,d^\bullet)$ is a complex of R-modules, then trivially
$(C^\bullet,\partial^\bullet_+=d^\bullet,\partial^\bullet_-=0)$ is a parity
complex. The $1^{st}$-order cocycles are the usual cocycles,
the $2^{nd}$-order cocycles are the cochains and the zero cobordisms
are the coboundaries.

As another example, the total complex $C^n=\underset{p+q=n}{\oplus} C^{p,q}$
of a double complex $(C^{p,q},d_1,d_2)$ may be viewed as a parity
complex with $\partial^+=d_1$ and $\partial^-=d_2$. Since $d_i^2=0$
the condition for the parity quasicomplex to be a complex $\partial^2=0$
is equivalent with $\{\partial^+,\partial^-\}=0$, the
double complex condition.
\begin{rem}\label{R:action}
If $(C^n,\mu_n,\partial^\pm_n)$ is a parity quasicomplex of monoids,
one may interpret the
maps $\partial^\pm_{n-1}$ as defining a biaction structure on $C^n$.
For $f\in C^{n-1}$ and $c\in C^n$
$$f\cdot c=(\delta^+f)\circ c,\qquad c\cdot f=c\circ(\delta^-f)$$
defines a left and a right action, with $\mu_n\circ_1\partial^+_{n-1}$ and
$\mu_n\circ_2\partial^-_{n-1}$ the corresponding maps.

Two elements are {\em cohomologous} if $f\cdot c=c'\cdot f$.
A {\em cocycle} is an element $f$ such that $f\cdot 1=1\cdot f$.
Two elements $c=f\cdot 1$ and $c'=1\cdot f$ are {\em cobordant}.
\end{rem}
%

\begin{example}\label{Ex:Hoch}
The additive Hochschild parity complex of an algebra $(A,\mu)$
is obviously equivalent to the Hochschild complex, with
$d=\partial^+-\partial^-$.
\end{example}

\subsection{Multiplicative Parity Quasicomplex}
In non-abelian cohomology one considers 
$G$ and $H$ two non-abelian groups and $\cdot:G\times H\to H$
an action. The cocycle and cohomology relations are defined
on cochains $C^n(G,H)$, which are normalized functions
$f:G\times...\times G\to H$ (section \ref{S:nacoh}).
Considering the group algebras $kG$ and $kH$,
over some commutative ring
$k$, the cohomology can be stated in terms of a
standard parity quasicomplex
of algebras analogous to the additive Hochschild parity quasicomplex.

Let $(A,\mu)$ and $(B,\cdot)$ be two arbitrary algebras (example
$A=kG$ and $B=kH$). Let $L:A\times B\to B$ and $R:B\times A\to A$ be two
$k$-linear functions. $B$ is referred to as an {\em $A$@-quasibimodule}.
With $f\in C^p(A;B)=Hom_k(A^{p+1},B)$\ ($p\ge-1$), define:
\begin{equation}
\partial^0_{p+2} f=L\circ_2 f,\qquad
\partial^i_{p+2} f=f\circ_i\mu,\qquad
\partial^{p+2}_{p+2} f=R\circ_1 f
\end{equation}
with less emphasis on the composition algebra structure. Define
\begin{align}
\partial^+=\prod_{i\ even}^{\to} \partial^i_{p+2}\qquad
\partial^-=\prod_{i\ odd}^{\leftarrow} \partial^i_{p+2}
\end{align}
\begin{defin}
In the context described above,
$(C^\bullet(A;B),\partial^\pm)$ is the
{\em standard multiplicative (Hochschild) parity quasicomplex
associated to the $A$-quasibimodule $B$}.
\end{defin}
For small values of $p$ we have:
\begin{alignat}{3}\label{Eq5}
&p=-1 \quad &\partial^+f&=L\circ_2 f &\partial^-f&=R\circ_1 f\notag\\
&\phantom{p=-1} \quad &( \partial^+&=\partial^0_1)
            &(\partial^-_1&=\partial^1_1)\notag\\
&p=0 \qquad &\partial^+f&=L\circ_2 f \cdot R\circ_1 f\qquad
            &\partial^-f&=f\circ_1\mu\notag\\
&\phantom{p=0} \qquad &( \partial^+&=\partial^0_2 \cdot \partial^2_2 )\qquad
            &( \partial^-&=\partial^1_2 )\notag\\
&p=1 \qquad &\partial^+f&=L\circ_2 f \cdot f\circ_2\mu \qquad
            &\partial^-f&=R\circ_1 f \cdot f\circ_1\mu \notag\\
&\phantom{p=1} \qquad &( \partial^+&=\partial^0_3 \cdot \partial^2_3 )\qquad
            &( \partial^-&=\partial^3_3 \cdot \partial^1_3 )\\
&p=2 \qquad &\partial^+f&=L\circ_2 f \cdot f\circ_2\mu
                                       \cdot R\circ_1 f\qquad
            &\partial^-f&=f\circ_3\mu \cdot f\circ_1\mu\notag\\
&\phantom{p=2} \qquad &( \partial^+
                    &=\partial^0_4 \cdot \partial^2_4
                                   \cdot \partial^4_4 )\qquad
            &( \partial^-&=\partial^3_4 \cdot \partial^1_4 )\notag
\end{alignat}
where $C^{-1}(A;B)=B$ ($A^0=\{0,1\}$) and $R\circ_1 b(a)=R(b,a)$ for
$p=-1$.



\section{\bf Cohomology of Monoidal Categories}\label{S:cohmc}
There is an analog construction for a functor
$S:(\C{C},\otimes)\to (\C{D},\otimes)$ between two categories
with products (not assumed monoidal).
The cochains are natural transformations.

The functor $F$ defines a {\em $\C{C}$@-quasibimodule} structure on $\C{D}$:
$$L=\tens{D}\circ_1F:\C{C}\boxtimes\C{D}\to\C{D},\quad
  R=\tens{D}\circ_2F:\C{D}\boxtimes\C{C}\to\C{D}$$
No constraints as to $L$ and $R$ being ``actions'' are assumed.
  
If $\phi\in C^p$ and $i=1,...,p+1$ define:
\begin{align}\label{E:delt}
(\delta^i_{p+2}\phi)_{A_1,...,A_{p+2}}&=
(\phi\circ_i \tens{C})_{A_1,...,A_{p+2}}=
  \phi_{A1,...,A_i\tens{C}A_{i+1},...,A_{p+2}},
  \notag\\
(\delta^{p+2}_{p+2}\phi)_{A_1,...,A_{p+2}}&=
(R\circ_1\phi)_{A_1,...,A_{p+2}}
   =R(\phi_{A1,...,A_{p+1}}, I_{A_{p+2}})\notag\\
   &=\phi_{A1,...,A_{p+1}}\tens{D} I_{F(A_{p+2})}\notag\\
(\delta^0_{p+2}\phi)_{A_1,...,A_{p+2}}&=
(L\circ_2\phi)_{A_1,...,A_{p+2}}
   =L(I_{A_{p+2}},\phi_{A1,...,A_{p+1}})\notag\\
   &=I_{F(A_{p+2})}\tens{D}\phi_{A1,...,A_{p+1}}\notag\\
\end{align}
Now for $p=0,1,2$, define:
\begin{align}
\partial^+\phi&=(\prod_{i\ even}^{\to} \partial^i_{p+2})\phi =
(\partial^0_{p+2}\phi)\circ(\partial^2_{p+2}\phi)
\circ(\partial^4_{p+2}\phi)\circ\dots\notag\\
\partial^-\phi&=
(\prod_{i\ odd}^{\leftarrow} \partial^i_{p+2})\phi
\dots\circ(\partial^3_{p+2}\phi)\circ(\partial^1_{p+2}\phi)
\end{align}
We interpret some categorical conditions in terms of parity
quasi-complexes.

%
Let $F:(\C{C},\otimes)\to (\C{D},\otimes)$ be a functor.
Recall that a {\em monoidal structure} is a 2-cochain:
$$\varPhi_{X,Y}:F(X\otimes Y)\to F(X)\otimes F(Y) \qquad
\text{for any} X,Y\in\C{C}$$
for which the following diagram is commutative for any objects
$X,Y,Z\in\C{C}$
\xycompileto{dia_p1}|{
\begin{gather}\diagram
\sst{F((X\otimes Y)\otimes Z)} \rto^{\varPhi_{X\otimes Y,Z}}
                         \dto_{F(\alpha_{X,Y,Z})}
& \sst{F(X\otimes Y)\otimes F(Z)} \rto^{\varPhi_{X,Y}\otimes I_{F(Z)}\qquad}
& \quad \sst{(F(X)\otimes F(Y))\otimes F(Z)} \dto^{\alpha_{F(X),F(Y),F(Z)}}\\
\sst{F(X\otimes (Y\otimes Z))} \rto^{\varPhi_{X,Y\otimes Z}} 
& \sst{F(X)\otimes F(Y\otimes Z)} \rto^{I_{F(X)}\otimes \varPhi_{Y,Z}\qquad}
& \sst{\quad F(X)\otimes (F(Y)\otimes F(Z))}\\
\enddiagram\label{D:mo1}\end{gather}}
More rigorously, $\Phi$ is a functorial morphism between the
following two functors:
$$\Phi:\partial^-F\to \partial^+F$$
(where the analogous maps $\partial^\pm$ on functors may be defined),
such that $\alpha$ and $\alpha'$ are $F$@-cohomologous
(\cite{I}, section 6):
\begin{gather}\diagram
F((X\tens{C}Y)\tens{C}Z)\quad
  \dto_{F(\alpha)} \rto^{\delta^-_3\Phi}
  & \quad (F(X)\tens{D}F(Y))\tens{D}F(Z) \dto^{\alpha'_F}\\
F(X\tens{C}(Y\tens{C}Z)) \quad \rto^{\delta^+_3\Phi}
  & \quad F(X)\tens{D}(F(Y)\tens{D}F(Z))\\
\enddiagram\label{D:mo3}\end{gather}
If the monoidal categories are strict, then the monoidal
structures of the functor $F$ are the 2-cocycles:
$\partial^+\Phi=\partial^-\Phi$.

\vspace{.1in}
Recall that a {\em monoidal morphism}
$\eta:(F,\varPhi)\to(G,\varGamma)$ of
monoidal functors $F,G:\C{C}\to\C{D}$ is a functorial morphism such that
the square diagram is commutative for any$X,Y\in\C{C}$:
\xycompileto{dia_p1}|{
\begin{gather}\diagram
F(X \otimes Y) \rto^{\varPhi_{X,Y}} \dto_{\eta_{X \otimes Y}}
& F(X) \otimes F(Y) \dto_{\eta_X \otimes \eta_Y}
& & \one_\C{D} \dlto_{\phi} \drto^{\gamma}\\
G(X \otimes Y) \rto^{\varGamma_{X,Y}} & G(X) \otimes G(Y)
& F(\one_\C{C}) \rrto^{\eta_{\one_\C{C}}} & & G(\one_\C{C})
\enddiagram\label{D:mo2}\end{gather}}
Equivalently, the first diagram is:
$$\diagram
\partial^-F \dto^{\partial^-\eta} \rto^{\Phi} &
  \partial^+F \dto^{\partial^+\eta}\\
\partial^-G \rto^{\Gamma} & \partial^+G
\enddiagram$$
i.e. $\Phi$ and $\Gamma$ are {\bf cohomologous} monoidal structures.
$$\diagram
\partial^+\eta \cdot \Phi = \Gamma \cdot \partial^-\eta
& \text{or} & \Gamma \overset{\eta}{\sim} \Phi
\enddiagram$$

\subsection{Cohomology of Monoidal Categories}\label{SS:cohmc}
We assume now that $F$ is the identity functor on the category
$\C{C}$ with product $\otimes$.

For each fixed {\em quasi-associator} $\alpha$, not assumed to be coherent,
we have the following parity quasicomplex (\cite{I}, definition 6.1):
\xycompileto{dia1}|{
\begin{gather}\diagram
 & & End(\otimes^2) \drto^{R_\alpha} \\
End(I) \ar@/^/[r]^{\partial^+} \ar@/_/[r]_{\partial^-} & End(\otimes)
\urto^{\partial^+} \drto_{\partial^-} & &
Hom(^2\otimes,\otimes^2) \ar@/^/[r]^{\partial^+} \ar@/_/[r]_{\partial^-}
 & Hom(^3\otimes,\otimes^3)\\
 & & End(^2\otimes) \urto_{L_\alpha} \\
C^0 \ar@/^/[r]^{d^0_+=\partial^+} \ar@/_/[r]_{d^0_-=\partial^-} &
C^1 \ar@/^/[rr]^{d^1_+=R_\alpha\circ\partial^+}
   \ar@/_/[rr]_{d^1_-=L_\alpha\circ\partial^-} & &
C^2 \ar@/^/[r]^{d^2_+=\partial^+} \ar@/_/[r]_{d^2_-=\partial^-} &
C^3
\enddiagram\label{D:monpqc}\end{gather}}
The notations $d^+_\alpha=d^1_+$ and $d^-_\alpha$ will be used to
stress the dependency on $\alpha$.

The {\em reduced parity quasi-complex} $(U^\bullet,d^\bullet_\pm)$
consists only of functorial isomorphisms.

Then we have (\cite{I}, corollary 6.2):
\begin{cor}
The 3-cocycles of the reduced parity quasicomplex
$(U^\bullet,d^\bullet_\pm)$ are the monoidal structures
of the category $\C{C}$ with product $\otimes$.
\end{cor}

\subsection{Commutativity Constraints}\label{SS:cc}
Specialize the category $\C{D}$ to the opposite category
$\C{C}_{op}$ (\cite{I}, section 6.4), and assume that their associators
are cohomologous.
Then any isomorphism $c$ conjugating the associator $\alpha$ and
its opposite $\alpha^{op}=c\alpha c^{-1}$ is called 
a {\bf commutativity constraint}.

We summarize the results on coboundary Hopf algebras (\cite{I}).
\begin{th}\label{T:cc}
Let $(H,\C{R})$ an almost cocommutative Hopf algebra and $\C{C}=H-mod$
the category of its representations. Then  the following statements
are equivalent\\
(i) $(H,\C{R})$ is a coboundary Hopf algebra.\\
(ii) $(H-mod, \sigma_\C{R})$ is a commutative monoidal category.\\
(iii) $(I,\sigma_\C{R}):H-mod\to H-mod_{op}$ is a monoidal equivalence.\\
(iv) $\sigma_\C{R}$ is a 2-cocycle.
\end{th}


\section{\bf Categorification}\label{S:categ}
An {\em affine group} with structure group $G$ is a set $E$ and a function
$\partial:E\times E\to G$, called {\em affine structure}, verifying
$\partial(b,c)\partial(a,b)=\partial(a,c)$.

In other words,
if $E$ is given the trivial groupoid structure $\C{C}_E$ - with
objects elements of $E$ and unique maps between objects
$Hom(e_1,e_2)=\{(e_1,e_2)\}$ - and $G$ is categorified as usual $\C{C}_G$
- the groupoid with one object $\C{O}b(\C{C}_G)=\{G\}$ and
$Hom(G,G)=G$ with group multiplication as composition of morphisms -,
then an affine group $E$ with structure group $G$ and affine structure
$\partial$ is just a (constant) functor $\partial:E\to G$.

We note that the standard categorification described above -
$\C{S}et \overset{Std}{\To} \C{C}at$ given by $E\mapsto \C{C}_E$ and
$\C{G}roups \overset{Std}{\To} \C{C}at$ given by $G\mapsto \C{C}_G$ -
is functorial.

We will be interested in a categorification procedure suited for
group extensions.

\subsection{Base Categorification}\label{SS:bc}
To a group $G$ associate the standard
monoidal category $\C{B}_G$ with simple objects
$\C{O}b(\C{B}_G)=\{G\}$ and $Hom(a,b)$ empty, unless $a=b$, when
the only morphism is the identity morphism. The monoidal product is
group multiplication. In an obvious way it is
the canonical skeletal monoidal groupoid $\C{C}$ with Grothendieck monoid
$\pi_0(\C{C})=G$ - isomorphism classes of objects/ connected components
of the associated geometric realization.

To a morphism of groups $s:G\to N$ associate the functor $S=\C{B}_s$
which on objects is just the morphism $s$. Since the only morphisms of
$\C{B}_s$ are identity morphisms,
this determines the value of $S$ on morphisms, being obviously a functor.

\subsection{Fiber Categorification}\label{SS:fc}
A group $G$ is a tautological affine group $E=G$ with structure group $G$
and affine structure
$$\partial:E\times E\to G,\qquad \partial(a,b)=b^{-1}a$$
There is a unique left and right action of $G$ on $G$ - $L^G$ and $R^G$ -,
associated with the natural left and right action of $G$ on the affine
group $E$ - $L^E$ and $R^E$ -, which is compatible with
the affine structure, i.e. commuting with $\partial$.
Indeed:
$$\partial(L^E_c(a),L^E_c(b))=\partial(ca,cb)=cb (ca)^{-1}=c(ba^{-1})c^{-1}
,\quad L^G_c(\partial(a,b))=L^G_c(ba^{-1})$$
$$\partial(R^E_c(a),R^E_c(b))=\partial(ac,bc)=bc (ac)^{-1}=ba^{-1}
,\quad R^G_c(\partial(a,b))=R^G_c(ba^{-1})$$
so that the left action of $G$ on itself is conjugation and
the right action is trivial:
$$L^G_c(x)=c x c^{-1}\qquad R^G_c(x)=x$$
The corresponding action on pairs is still denoted as
$L_c(a,b)=(L_c(a),L_c(b))$ and the superscripts $G/E$ will be omitted.

Define a monoidal category $\C{F}_G$
with objects $\C{O}b(\C{F}_G)=G$ and morphisms $Hom(a,b)=G$,
elements of $G$.
The composition of morphisms is multiplication in $G$ and 
the monoidal product is defined as follows:
\begin{gather}
\diagram
a \dto^{f \quad\otimes} &
  b \dto^{g \qquad =\quad} & ab \dto^{t(f)gs(f)^{-1}=a' g a^{-1}}\\
a' &  b' & a' b' &
\enddiagram\label{E:mp1}\end{gather}
On objects the monoidal product is again multiplication in $G$,
but on morphisms, it is
``left twisted'': $f\otimes g=t(f) g s(f)^{-1}$. Here $t$ and $s$
denote the target and source maps. In particular, if
$f$ is the identity morphism $I_c$, then:
$$\diagram
c \dto^{1 \quad\otimes} &
   a \dto^{g \qquad =\quad} & ca \dto^{c g c^{-1}}\\
c & b & c b &
\enddiagram$$
and if $g$ is the identity morphism $g=I_c$, then:
$$\diagram
a \dto^{f \quad\otimes} &
  c \dto^{1 \qquad =\quad} & ac \dto^{b a^{-1}}\\
b &  c & bc &
\enddiagram$$
It is easy to see that $\otimes$ is a functor:
\begin{align}
\otimes((g,g')\circ(f,f'))&=\otimes((g\circ f, g'\circ f'))
=t(g)(g'f')s(f)^{-1}\notag\\
\otimes(g,g')\circ\otimes(f,f')&=t(g)g's(g)^{-1} t(f)f's(f)^{-1}
=t(g)g's(f)^{-1}\notag\\
\otimes((I_a,I_b))&=a 1 a^{-1}=1\notag
\end{align}
In the above equations we have omitted to show the obvious
source and target of the various maps involved.
\newcommand{\ve}[1]{\overset{\to}{#1}}
The twisted monoidal product corresponds to the action of $G$ on itself,
when interpreted as an affine space. 

The morphism $a\overset{ba^{-1}}{\to}b$ in $Hom(a,b)$ is called
the {\em vector} $\ve{(a,b)}$ (or $\ve{ab}$ if no confusion is possible)
and $\ve{ }$ is a cross-section of the fibration
$Hom \overset{(s,t)}{\To} \C{O}b\times \C{O}b$, also denoted by $\partial$.
Then, the right translation
$\ve{ab}\cdot c=\ve{(ac,bc)}$ changes only the source and the target, since
$(bc)(ac)^{-1}=ba^{-1}$,
and corresponds to the monoidal multiplication
to the right with the identity morphism $I_c$:
$$\diagram
a \dto_{f=ba^{-1}}^{\quad\otimes} &
   c \dto^{1 \quad =\quad} & ac \dto^{b a^{-1}=f}\\
b & c & bc &
\enddiagram$$
The left translation conjugates the vector:
$c\cdot \ve{ab}=\ve{(ca,cb)}=cba^{-1} c^{-1}=I_c\cdot\ve{ab}$,
and corresponds to the monoidal multiplication
to the left with the identity morphism $I_c$:
$$\diagram
c \dto_{1}^{\quad\otimes} &
   a \dto^{ba^{-1} \ =\quad} & ca \dto^{c (ba^{-1}) c^{-1}}\\
c & b & c b &
\enddiagram$$
Note that right multiplication of an arbitrary morphism
$a\overset{f}{\to}b$ with
the identity morphism $I_c$ translates the objects and
truncates (projects) the morphism onto the corresponding vector
$a\overset{ba^{-1}}{\to}b$.

\begin{rem}\label{R:1}
The cross-section $\partial$ intertwines the left
and right action of $G$ on $Hom$ given by left and right
monoidal multiplication with identity morphisms,
with the natural actions of $G$ on $E=G$ as an affine group:
$$I_c\otimes \partial(a,b)=\partial(c\cdot(a,b)),\qquad
(I_c\otimes \ve{(a,b)}=\ve{c(a,b)})$$
$$\partial(a,b)\otimes I_c=\partial((a,b)\cdot c),\qquad
(\ve{(a,b)}\otimes I_c=\ve{(a,b)c})$$
Also note that $\partial$ is functorial:
$$\partial((b,c)\circ(a,b))=\partial(b,c)\circ\partial(a,b)$$
\end{rem}
For reference, we state:
\begin{cor}\label{cor:tens}
Any categorical diagram in $\C{F}_G$ with vectors as morphisms commutes.

Left or right tensoring with identity maps is vector preserving.
\end{cor}
Note that all objects are isomorphic, so that the Grothendieck monoid
$\pi_0(\C{F}_G)$ is trivial.

\vspace{.2in}
{\bf The functor $\C{F}$.}
We will define $\C{F}$ on morphisms.

The goal is
to extend $\C{F}$ to the enlarged category $\tilde{\C{G}r}$
with objects groups, but having {\bf functions as morphisms}.
Assuming $s:G\to N$ is an arbitrary function, then $\C{B}$ should
associate to it a monoidal functor with monoidal structure
$s(ab)\overset{f(a,b)}{\to}s(a)s(b)$, with $f(a,b)=s(a)s(b)s(ab)^{-1}$.
Since $\C{F}_G$ is a groupoid, the coboundary map
$\delta(s)=\partial^+s \cdot (\partial^-s)^{-1}$ may be defined for
functors (section \ref{S:cohmc}). Then $f=\delta s$.
In section \ref{S:nacoh} will be shown that it equals
the group cohomological coboundary $f=\delta_L s$,
where $L=C_s:G\to Aut(N)$ is the quasi-action of $G$ on $N$ defined 
through the conjugation induced by $s$.

%
\newcommand{\ff}{$\C{F}'$}
\newcommand{\pair}{$G\underset{\Lambda}{\overset{s}{\rightrightarrows}}N$}
We will first define for each pair \pair of functions
a functor $F(s,\Lambda):\C{F}_G\to \C{F}_N$ by
$F(s,\Lambda)(x:a\to b)=(z:s(a)\to s(b))$,
where $z=s(b) \Lambda(b^{-1}xa) s(a)^{-1}$.
\begin{lem}
$F(s,\Lambda)$ is a functor iff $\Lambda$ is a morphism of groups.
\end{lem}
We have two natural definitions of $\C{F}$ on the category $\C{G}r$ of
groups, corresponding to a constant second component $(s,1)$ or to the
diagonal of $\C{G}r\times \C{G}r$:
$$\C{F}(s)=F(s,1), \qquad \C{F}_{\Delta}=F(s,s)$$
\begin{defin}
Define $\tilde{\C{G}}r$ the extended category of groups, with objects
groups and morphisms functions preserving the group identities.
\end{defin}
\begin{prop}\label{P:functorf}
The functor $\C{F}$ naturally extends to $\tilde{\C{G}}r$:
$$ \C{F}(s)=s^*(\partial_N), \quad s\in Hom_{Set}(G,N)$$
Moreover $f=\delta s$ is a monoidal structure for $\C{F}(s)$, so that
$\C{F}$ is valued in the category $\C{M}on$ of monoidal categories.

$\C{F}(s)$ is a strict monoidal functor iff $s$ is a group morphism.
\end{prop}
\begin{pf}
Denote by $(S,s)=\C{F}(s)$, with $S$ the functor component acting on
functions. Recall that
$S(x:a\to b)=(z:s(a)\to s(b))$ with $z=s(b)s(a)^{-1}$,
and
$(s^*\partial_N)(a,b)=\partial_N(s(a),s(b))=s(b)s(a)^{-1}:s(a)\to s(b)$.

That $f(a,b)=\delta_Ls(a,b)=s(a)s(b)s(ab)^{-1}$ is a functorial morphism:
$$\diagram
s(ab) \dto_{S(x\otimes y)} \rto^{f(a,b)} & s(a)s(b) \dto^{S(x)\otimes S(y)}\\
s(a'b') \rto_{f(a',b')} & s(a') s(b') 
\enddiagram$$
where $x:a\to b$ and $y:a'\to b'$ are two arbitrary morphisms
in $\C{F}_G$, follows by a direct computation.

That $f$ is a monoidal structure:
$$\diagram
s((xy)z)\quad \rto^{f(xy,z)} \dto^{s(1)} &
  \quad s(xy)s(z)\quad \rto^{f(x,y)\otimes I_{s(z)}}
  & \quad (s(x)s(y))s(z) \dto^{1}\\
s(x(yz))\quad \rto^{f(x,yz)} & \quad s(x)s(yz)\quad
  \rto^{I_{s(x)}\otimes f(x,y)} & \quad s(x)(s(y)s(z))
\enddiagram$$
follows, since $f(x,y)$ is the vector $s(xy)\to s(x)s(y)$ and
left or right tensoring with $I_c$ is vector preserving
(corollary \ref{cor:tens}).
\end{pf}
This proves the functoriality of the above construction. 

\vspace{.2in}
{\bf Untwisting $\C{F}_G$.}
To better understand fiber categorification, consider the
underlying category $\C{F}_G$ with monoidal product $\tilde{\otimes}$
defined as group multiplication on objects (same as $\otimes$) and as
the projection on the second factor on morphisms:
$$\diagram
a \dto^{f \quad\otimes} &
  b \dto^{g \qquad =\quad} & ab \dto^{g}\\
a' &  b' & a' b' &
\enddiagram$$
Define $D:(\C{F}_G, \otimes)\to (\C{F}_G, \tilde{\otimes})$ as identity on
objects and associating to a morphism $a\overset{x}{\to}b$ the ``loop
based at 1'':
$$\diagram
D(a\overset{x}{\to}b)=a\overset{\tilde{x}}{\to}b
\ar@{}[d]|{\tilde{x}=b^{-1}xa} & a \rto^{x} & b\\
 & 1 \uto_{a} \urto_{b}
\enddiagram$$
\begin{prop}\label{P:fiber}
$(\C{F}_G, \tilde{\otimes})$ is a strict monoidal category and
$D$ is a strict monoidal isomorphism.
\end{prop}
\begin{pf}
On morphisms $\tilde{\otimes}$ is defined as follows:
$$\diagram
a \dto^{\tilde{x} \quad\otimes} &
  b \dto^{\tilde{y} \qquad =\quad} & ab \dto^{\tilde{y}}\\
a' &  b' & a' b' &
\enddiagram$$
where $\tilde{x}:a\to a'$ and $\tilde{y}:b\to b'$ are morphisms in
$(\C{F}_G, \tilde{\otimes})$.

It is easy to see that $\tilde{\otimes}$ and $D$ are functors. 

Now only note that $D(x\otimes y)=D(x)\tilde{\otimes}D(y)$. To see this,
$D(x\otimes y)=D(a'ya^{-1}:ab\to a'b')=\tilde{z}:ab\to a'b'$, where
$\tilde{z}=(a'b')^{-1} a'ya^{-1} (ab)=(b')^{-1}ya'=\tilde{y}$.
But $D(x)\tilde{\otimes} D(y)=\tilde{y}:ab\to a'b'$.
\end{pf}

\subsection{Bundle Categorification}\label{SS:bdlcat}
Let $\C{E}$ be an extension $1\to N\to E\to G\to 1$.
We associate  a category $\C{C}(\C{E})$ in the following way.
Think of $E$ as principal bundle over $G$ with structure group $N$.
Use fiber categorification for its affine fibers.
Morphisms are $f:a\to b$ with $a,b\in E_g$, the fiber over $g\in G$,
and $f\in N$. The monoidal product is defined as in \ref{E:mp1}:
$$\diagram
a \dto^{f \quad\otimes} &
   b \dto^{g \qquad =\quad} & ab \dto^{t(f)gs(f)^{-1}=a' g a^{-1}}\\
a' & b' & a' b' &
\enddiagram$$
Let $s:G\to E$ be a set theoretic section $s:G\to E$
(all are assumed normalized: $s(1)=1$, i.e. a morphism in \Gr),
trivializing $E$ as a set.
If $a=s(x)n$ and $a'=s(x)n'$ with $x\in G$, then
$a' g a^{-1}=s(x)n' x n^{-1}s(x)^{-1}$ is still an element of $N$.
Note that the monoidal product respects fibers, since if $g$ is
over $y\in G$ ($g\in Hom(a',b')$ with $a', b'\in E_y$ )
then $f\otimes g$ is a morphism over $xy$.
\begin{rem}\label{R:vect}
More importantly, the monoidal product preserves vectors:
$$\partial(a,a')\otimes \partial(b,b')=
\partial(ab, a'b')$$
If $f=a'a^{-1}$ and $g=b'b^{-1}$ then $a'ga^{-1}=a'b'b^{-1}a^{-1}$
and the morphism in the right hand side of the above diagram is the
corresponding vector.
\end{rem}
\begin{th}\label{T:categ}
Base, Fiber and Bundle Categorification maps $\C{B}, \C{F}$ and
$\C{C}$ are functors from
the category of groups $\C{G}r$ and extensions $\C{E}$ to the category
$\C{M}on_s$ of strict monoidal categories.

$\C{B}$ and $\C{F}$ are the restrictions of $\C{C}$
corresponding to the two natural functors:
$$ G \overset{T_b}{\mapsto} (1\to G\to G\to 1\to 1)\quad (trivial\ base)$$
$$ G \overset{T_f}{\mapsto} (1\to 1\to G\to G\to 1)\quad (trivial\ fiber)$$
embedding the category of groups into the category of group extensions:
$$
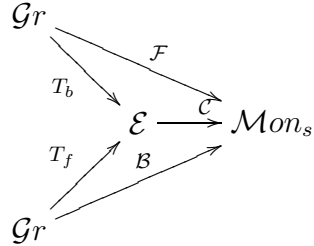
\diagram
\C{G}r \drto_{T_b} \drrto^{\C{F}} \\
& \C{E} \rto^{\C{C}} & \C{M}on_s \\
\C{G}r \urto^{T_f} \urrto^{\C{B}}
\enddiagram$$
\end{th}

If we disregard the change of coeficients in a categorification,
the base categorification corresponds to $K-categorification$
in the sense of \cite{CY2}.

\begin{rem}
According to \cite{Porter}, example page 692, ``The concept of a
$cat^1$@-group is equivalent to that of a crossed module.''
Recall that a {\em crossed module} is
an $E$-equivariant morphism $\phi:N\to E$ of $E-groups$ - with inner
conjugation of $E$ on $E$ - such that the ``composition'' of
the left action of $N$ on $E$ through $\phi$ with the $E$ action on
$N$ $L:N\to Aut(E)$, is the inner conjugation $C:N\to Aut(N)$:
$(n\cdot 1_E)\cdot m=n\cdot(1_E\cdot m)$. Or the following diagram
commutes (see \cite{Brown}):
$$\diagram
N \drto_{C} \rto^{\phi} & E \dto^{L}\\
& Aut(N) \\
\enddiagram$$
Since a crossed module defines a 4-term exact sequence
$$ 0\to A\to N\overset{\phi}{\to} E\to G\to 1$$
with $A=ker\phi$ and $G=coker\phi$ (see \cite{Brown}, section 4.5, page 102),
one expects a definite relation between the categorification of an
extension as described above, and $\C{C}at^1$@-groups (see also
\cite{Breen}). It will not be discussed herein.
\end{rem}


\section{\bf Non-abelian Cohomology of Groups}\label{S:nacoh}

Let $G$ and $N$ be groups.

When defining cohomology with non-abelian coefficients $N$ we
consider all $G$@-structures on the group $N$ and relax
the action requirement. In other words, we work in the extended group
category \Gr with groups as objects, but functions respecting the
identity element as morphisms.
The categorical interpretation of functions between groups,
made precise in section \ref{S:categ}, is that of monoidal functors.
Such a monoidal functor is strict if the corresponding function
is a 1-cocycle, i.e. it is a group morphism. When a 
confusion is possible, we will call a $\C{G}r$@-morphism a {\em strict
morphism} when viewed in \Gr (and a \Gr@-morphism,  a function!).

The cartesian power $G\times ...\times G$ of $G$ $n$@-times, is denoted as
$G^n$. Let
$$C^k(G,N)=\{(f,L)| f:G^k\to N, \ L:G\to Aut(N),
  \quad f, L \ normalized \ functions\}.$$
be the set of normalized k-cochains, i.e. the natural "multi-morphism" maps
in the category \Gr: $f(...,1,...)=1$ and $L(1)=1$. In what follows we will
use the unshifted grading.

$L:G\to Aut(N)$ - a normalized function - is a morphism in \Gr,
so it is an action in \Gr. To avoid possible confusions, it will be called
a {\em quasiaction}.

We think of $\coprod_kG^k$ as corresponding to the bar resolution in
\Gr and of $C^\bullet(G,N)=\coprod_kC^k(G,N)$ as the collection of all
corresponding (multiplicative) Hochschild cochains.
Here ``all'' means, for all possible $G-group$ structures on $N$. Recall
that a $G$-group structure $L:G\to Aut(N)$ in \Gr is just a function, and
need not be a morphisms of groups.

Define the functions $\delta^k:C^k\to C^{k+1}$ by:
$$ \delta^k(f,L)=(\delta^k_L(f),L)$$
where $\delta^k_L$ is the usual coboundary map of the multiplicative
parity quasicomplex $C^\bullet(G, _LN)$, as defined in section
\ref{S:pqc}. It corresponds to the left \Gr-action of $G$ on $N$ defined
by $L$, and to the right trivial action.
We recall the formulas and use the unshifted degrees.
The dependence on $L$, the {\em $G$-quasiaction} on $N$, is included since
the structure introduced later, will involve several {\em $L$-fibers}
$C^\bullet(G, _LN)$ (the cochains with a fixed $L$).

If $(f,L)\in C^p(G,N)$:
$$
\partial^0_{p+1} (f,L)=(L\circ_2 f,L)\qquad
\partial^i_{p+1} (f,L)=(f\circ_i\mu,L)\qquad
\partial^{p+1}_{p+1} (f,L)=(R\circ_1 f,L)
$$
\begin{align}\label{E:pm}
\partial^+=\prod_{i\ even}^{\to} \partial^i_{p+1}\qquad
\partial^-=\prod_{i\ odd}^{\leftarrow} \partial^i_{p+1}\\
\delta=\partial^+\cdot (\partial^-)^{-1}\qquad
\bar{\delta}=(\partial^-)^{-1}\cdot \partial^+\notag
\end{align}
To abbreviate the notation we use:
$$\partial^\pm(f,L)=(\partial^\pm_L f,L),\qquad
\delta(f,L)=(\delta_Lf,L),\qquad \bar{\delta}(f,L)=(\bar{\delta}_Lf,L)$$
Explicitly, for $n\in C^0(G,N)=N,\ s\in C^1(G,N),\ f\in C^2(G,N)$
and $\alpha\in C^3(G,N)$, we have:
\begin{alignat}{3}\label{Eq10}
&p=0\quad & &\partial^+n(a)=L_a(n)\qquad \partial^-n(a)=n\notag\\
&p=1\quad & &\partial^+s(a,b)=L_a(f(b))\qquad \partial^-s(a,b)=s(ab)\notag\\
&p=2\quad & &\partial^+f(a,b,c)=L_a(f(b,c)) f(a,bc)\\
&         & &\partial^-f(a,b,c)=f(a,b)f(ab,c)\notag\\
&p=3\quad & &\partial^+\alpha(a,b,c,d)=L_a(\alpha(b,c,d))
  \alpha(a,bc,d)\alpha(a,b,c)\notag\\
&         & &\partial^-\alpha(a,b,c,d)=\alpha(a,b,cd)\alpha(ab,c,d)\notag
\end{alignat}
Note that $\delta$'s and $\partial$'s are \Gr@-morphisms (functions)
in a tautological way.
\begin{defin}
$(C^\bullet(G,N),\partial^\pm)$ is the standard parity quasicomplex
of $G$ with coefficients in $N$.

The set of cocycles $(f,L)$, i.e. verifying $\partial^+_Lf=\partial^-_Lf$,
is denoted by $Z^k(G,N)$.
The set of coboundaries $(\delta_Lf,L)$ is denoted by $B^k(G,N)$.
\end{defin}
\begin{notation}
The strict morphism defined by inner conjugation will be denoted with
$C:N\to Aut(N)$. By abuse of notation $C$ will also denote
$C:Aut(N)\to Aut(Aut(N))$.

Also the notations $C_n(x)=n x n^{-1}$ and $L_g(n)=L(g,n)$
will be used. $L_a^{-1}$ and $C_a^{-1}$ should be understood
as inverses in $Aut(N)$, i.e. $(L_a)^{-1}$.

If $(s,L)\in C^1(N,G)$ is a 1-cochain, then $f$ will usually denote
$f=\delta_Ls$ and referred to as the coboundary of $(s,L)$, the pair
$(f,L)$ being understood.

The $Aut(N)$ action $C\circ L$ of $G$ on $Aut(N)$ induced by $L$,
will be denoted by $\C{L}$.
\end{notation}
The morphism $C$ induces a function
$C_*:C^\bullet(N)\to C^\bullet(G,Aut(N))$
on the corresponding parity quasicomplexes:
$$ C_*(f,L)=(C_f, C_L), \quad C_f(g)=C_{f(g)}, \quad C_L(g)=C_{L(g)}$$
where $(f,L)\in C^k(G,N)$.
The grading notation for this map will be omitted.
%
%
\begin{lem}\label{L:chain}
$C_*:C^\bullet(G,N)\to C^\bullet(G,Aut(N))$ is a chain map.
\end{lem}
\begin{pf}
Note that $C_*$ commutes with the coface maps $\partial^i_{p+1}$
for $i=1,...,p$:
$$C_*(\partial^if)=C\circ (f\circ_i \mu)=(C\circ f)\circ_i \mu).$$
To check $\partial^0(C\circ f,C\circ L)=C\circ(\partial^0(f,L))$,
or equivalently $C\circ \partial^0_L(f)=\partial^0_\C{L}(C\circ f)$,
with $\C{L}=C_L$, note that
$$\C{L}_a(C(n))=C_{L_a}(C(n))=L_a\circ C_n\circ L_a^{-1}=C_{L_a(n)}
=C(L_a(n))$$
Also since $C$ is a morphism of groups $C\circ (\prod \partial^i)=
\prod C\circ\partial^i$, and $C_*$ commutes with $\partial^\pm$.
\end{pf}
For any group $N$ there is an associated exact sequence:
\begin{equation}\label{D:cen}
0\to Cen(N) \overset{i}{\to} N\overset{C}{\to} Aut(N)
\overset{\pi}{\to} Out(N) \to 1
\end{equation}
where  $Out(N)=Aut(N)/Int(N)$ is the quotient of $Aut(N)$
modulo inner actions.
\begin{prop}
$$0\to C^\bullet(G,Cen(N)) \overset{i}{\to}
C^\bullet(G,N)\overset{C_*}{\to} C^\bullet(G,Aut(N))
\overset{\pi}{\to} C^\bullet(G,Out(N)) \to 1$$
is an exact sequence of parity quasicomplexes.
\end{prop}
\begin{pf}
Since for $g\in G$ $L_g$ are morphisms of groups,
the inclusion $j$ and projection $\pi$ commute with the coboundary maps
$\partial^\pm$. $C_*$ is a chain map by the above lemma.
\end{pf}
The curvature $(\delta^\bullet)^2$ of the quasicomplex is investigated
in the following.
\begin{defin}
If $(c,L)$ is a p-cochain, define $\C{I}(c,L)$ ({\em holonomy group}
of $(c,L)$) as the
subgroup of $N$ generated by the orbit of $Im(c)$ under the quasiaction $L$:
$$ \C{I}(c,L)=< \{L_a(L_b(...(c(g_1,...,g_p))...)| a,b...,g_i\in G\}>.$$
An equation involving automorphisms of $N$ will be said to hold on $\C{I}$
(or on $c$) if it is verified when restricted to $\C{I}(c,L)$.
\end{defin}
\begin{lem}\label{L:holon}
Let $(c,L)\in C^p(N,G)$ be a p-cochain. Then
$\C{I}(c,L)$ is invariant under $L$, which naturally induces a quasiaction
$L^{|c}:G\to Aut(\C{I})$.

If $(s,L)$ is a 1-cochain and $L'=C_s^{-1}\circ L$, then $\C{I}(c,L)$
is also $L'$ invariant.
\end{lem}
\begin{lem}\label{L:dds}
For any 1-cochain $(s,L)\in C^1(G,N)$ and $f=\delta_Ls$,
the following conditions are equivalent:
\begin{alignat}{2}
&(i) \hspace{1in} & &\delta^2_L\circ \delta^1_L(s)=1\\
&(ii) & &\delta_\C{L}L=C_f \quad \text{on}\ \C{I}\\
&(iii) & &\delta_\C{L}L=\delta_\C{L}(C_s) \quad\text{on}\ \C{I}
  \label{E:iii}\\
&(iii) & &\gamma=C_s^{-1}\circ L \quad
  \text{on}\ \C{I},\ \text{is\ a\ morphisms\ of\ groups\quad}
\end{alignat}
where $\C{L}=C_L$.
\end{lem}
\begin{pf}
For $s:G\to N$ and $L:G\to Aut(N)$:
\begin{equation}\label{E:p1}
\partial^+_Ls(a,b)=L_a(s(b))s(a) \qquad \partial^-_Ls(a,b)=s(ab)
\end{equation}
\begin{equation}\label{E:p2}
f=\delta_L s=(\partial^+_Ls)(\partial^-_Ls)^{-1},\quad
f(a,b)=L_a(s(b))s(a)s(ab)^{-1}
\end{equation}
\begin{equation}\label{E:p3}
\partial^+_Lf(a,b,c)=L_a(f(b,c))f(a,bc),\quad
\partial^-_Lf(a,b,c)=f(a,b)f(ab,c)
\end{equation}
\begin{align}
\partial^+_Lf(a,b,c)
  &=L_a(L_b(s(c))s(b)s(bc)^{-1})L_a(s(bc))s(a)s(a(bc))^{-1}\notag\\
  &=L_a(L_b(s(c)))L_a(s(b))s(a)s(abc)^{-1}\notag\\
\partial^-_Lf(a,b,c)&=L_a(s(b))s(a)s(ab)^{-1}L_{ab}(s(c))s(ab)s((ab)c)^{-1}
\notag
\end{align}
Then $\partial^+_Lf=\partial^-_Lf$ iff:
\begin{align}
L_a (L_b(s(c)))L_a(s(b))s(a)&=L_a(s(b))s(a)s(ab)^{-1}L_{ab}(s(c))s(ab)
  \notag\\
L_a (L_b(s(c)))&=
  L_a(s(b))s(a)s(ab)^{-1}L_{ab}(s(c))s(ab)s(a)^{-1}L_a(s(b))^{-1}
\notag\\
&=C_{L_a(s(b))s(a)s(ab)^{-1}}(L_{ab}(s(c)))\notag\\
&=C_{f(a,b)}(L_{ab}(s(c)))\notag\\
L_a\circ L_b &= C_{\delta_L s}\circ L_{ab} \qquad
  \text{on\quad $\C{I}(s,L)$}
\end{align}
Let $\C{L}=C_L:G\to Aut(Aut(N))$ be the inner quasi-action
on $Aut(N)$, induced by $C$.
Note that, for the pair $(L,\C{L})\in C^1(G,Aut(N))$, we have:
$$\partial^+_\C{L} L(a,b)=L_a\circ L_b,\quad \partial^-_L L(a,b)=L_{ab}$$
$$\delta_\C{L} L(a,b)=(L_a\circ L_b) \circ L_{ab}^{-1}.$$
So $\delta^2_\C{L}\circ\delta^2_\C{L}=1$ iff
$\partial^+_\C{L} L=C_f \circ \partial^-_\C{L} L$
on $\C{I}$, i.e. $\delta_{C_L} L=C_{\delta_L s}$.

$C_*$ is a chain map, so that $C_*(\delta(L,s))=\delta(C_L,C_s)$.
Then $\delta^2_L\circ\delta^1_L=1$ iff $\delta(C_L,L)=\delta(C_L, C_s)$, or
$\delta_\C{L}(L)=\delta_\C{L}(C_s)$
\begin{gather}
\diagram
  (s,L) \ar@{|->}[d]^{C} \ar@{|->}[r]^{\delta^1} &
  (f=\delta_Ls,L) \ar@{|->}[d]^{C}\\
  (C_s, \C{L})  \ar@{|->}[r]^{\delta^1\quad} &
  (\delta_\C{L}(C_s),\C{L})=(C_f,\C{L})
\enddiagram\label{D:2}\end{gather}
Note that the condition \ref{E:iii} is equivalent to:
\begin{align}
L_a\circ L_b\circ L_{ab}^{-1}&=
  \C{L}_a\circ (C_s(b))\circ C_s(a)\circ C_s(ab)^{-1}
\quad \text{on}\ \C{I}\notag\\
&= L_a\circ C_{s(b)}\circ L_a^{-1}\circ C_{s(a)}\circ C_{s(ab)}^{-1}
\notag
\end{align}
and simplifies to:
$$L_b L_{ab}^{-1}=C_{s(b)} L_a^{-1} C_{s(a)} C_{s(ab)}^{-1}
\quad \text{on}\ \C{I}.$$
A direct computation:
$$C_{s(b)}^{-1}\circ L_b=L_a^{-1}\circ C_{s(a)}\circ C_{s(ab)}^{-1}
\circ L_{ab}$$
$$C_{s(a)}^{-1}\circ L_a\circ C_{s(b)}^{-1}\circ L_b=
C_{s(ab)}^{-1}\circ L_{ab}$$
shows it is equivalent to
$\gamma_a\circ \gamma_b=\gamma_{ab}$. 
\end{pf}
\begin{cor}\label{C:1}
If $(s,L)$ is a 1-cocycle then $L$ and $L'=C_s^{-1}\circ L$
corestrict to actions $L^{|s}, {L'}^{|s}:G\to Aut(\C{I})$.
\end{cor}
\begin{pf}
Since $\delta^2_L\circ \delta^1_L(s)=1$, the corestriction of $L'$ is a morphism
of groups by the above lemma. Or directly (see the remark below), since
$L_a(s(b))s(a)=\partial^+_Ls(a,b)=\partial^-_Ls(a,b)=s(ab)$:
\begin{align}
s(a(bc))&=L_a(s(bc))s(a)=L_a(L_b(s(c))s(b))s(a)=
  L_a(L_b(s(c))L_a(s(b))s(a)\notag\\
&=L_a(L_b(s(c)))s(ab)\notag
\end{align}
and
$$s((ab)c)=L_{ab}(s(c))s(ab).$$
\end{pf}
\begin{rem}\label{R:order}
Note that $\partial^+_Ls(a,b)=s(a)L'_a(s(b))s(a)^{-1}s(a)=s(a)L'_a(s(b))$.

It is convenient to {\em represent} a 1-cochain $(s,L)$ as $(s,C_s\circ L')$.
\end{rem}
\begin{lem}\label{L:p0}
Let $(n,L)$ be a 0-cochain. Then
$\delta^2_L\circ \delta^1_L\circ \delta^0_L(n)=1$ iff
$\delta_\C{L}L=C_f$ on $\C{I}$, where $s=\delta_Ln,\ f=\delta_Ls$
and $\C{I}$ is the subgroup generated by the orbit of $n$ under the
quasiaction $L$.

In particular, $\delta^1_L\circ \delta^0_L(n)=1$ iff $L$ restricts to an
action of $G$ on $\C{I}$.
\end{lem}



The equation \ref{E:MC}:
\begin{equation}
\partial^+_\C{L}L=C_f\circ \partial^-_\C{L}L \tag{MC}
\label{E:MC}\end{equation}
controls the curvature $\delta^2$ of the quasicomplex.
\begin{defin}\label{Def:int}
A p-cochain $(c,L)$ ($0\le p\le 2$) is {\em integrable} if equation
(MC) is verified on the corresponding subgroup $\C{I}=<L,c>$,
where $f=c$ if $p=2$, $f=\delta_Lc$ if $p=1$ and
$f=\delta^1_L\circ\delta^0_Lc$ if $p=0$.
It is called {\em absolute integrable} if the equation (MC)
holds globally on $N$.

A 1-cochain $(s,C_s)$ is called a {\em standard 1-cochain}.
\end{defin}
The set of standard 1-cochains is the graph of the morphism
$C_*:C^1(G,N)\to C^1(G,Aut(N))$. 
Note that if $(s,L=C_s)$ is a standard 1-cochain, then $(L,\C{L}=C_L)$
is also a standard cochain.
\begin{prop}\label{P:int}

(i) For $p=0,1$, p-cocycles are integrable.

(ii) A standard 1-cochain is absolute integrable.

(iii) If $(s,L)$ is a 1-cochain and $L'=C_s^{-1}\circ L$, then
$\C{I}(s,L)=\C{I}(s,L')$. The 1-cochain $(s,L)$ is (absolute) integrable
iff $L$ is a morphism of groups when corestricted to $\C{I}(s,L')$ ($N$).
\end{prop}
\begin{pf}
(i) The 0-cocycles $(n,L)$ are the fixed points $n$ of the quasiaction $L$:
$$L_g(n)=n,\quad any \ g\in G.$$
and $L$ restricted to $\C{I}(n,L)=<n>$ is a morphism of groups. Then
equation (MC) holds on $\C{I}$.

For a 1-cocycle $(s,L)$, $L^{|s}$ is an action by corollary \ref{C:1}.

(ii) Note that, since $C_*$ is a chain map,
$\delta_{C_L}L=C\circ \delta_Ls=C_f$ holds globally on $N$.

(iii) If equation (MC) holds on $\C{I}(s,L')=\C{I}(s,L)$,
then $\delta_Lf=1$ by
lemma \ref{L:dds}, where $f=\delta_Ls$ and $L=C_s\circ L'$. Then
$\delta_\C{L}L=1$ on $\C{I}$.

If $L'$ is a morphism on $\C{I}$ then a direct computation shows
$\partial^+_Lf=\partial^-_Lf$, and so $(f,L)$ is a cocycle. Then
by lemma \ref{L:dds} $\delta_\C{L}L=C_f$ holds on $\C{I}$, so that
$(s,L)$ is integrable.
\end{pf}
%
\begin{defin}
For $p=0,1,2$, $C^p_{int}(G,N)$ is the set of integrable cochains.
\end{defin}
A motivation for the terminology and an interpretation of
equation (MC) are given latter, using the categorical interpretation
(theorem \ref{T:MC}).

By proposition \ref{P:int}, we have:
\begin{cor}
$Z^0_{int}(G,N)=Z^0(G,N)$ and $Z^1_{int}(G,N)=Z^1(G,N)$.
\end{cor}
%

\subsection{Quasiaction of $C^k$ on $C^{k+1}$}\label{SS:action}

So far $C^\bullet(G,N)$ is just the disjoint union of the parity
quasicomplexes corresponding to individual quasiactions $L$.
Note that $(f_1,L_1)$ and $(f_2,L_2)$ are cohomologous (see
{\em parity complexes}) only if
$L_1=L_2$ and there is a one degree less cochain $(s,L)$ ($L=L_1$)
such that $\partial^+_Ls \cdot f=f'\cdot \partial^-_Ls$.
We will need a more general relation on cochains.

\begin{defin}
Two $p$@-cochains $(f,L)$ and $(f',L')$ are {\em weak cohomologous},
and denoted $(f,L)\underset{wk}{\sim}(f',L')$,
if there is a $(p-1)$@-cochain $(\gamma,l)$ such that:
$$\partial^+_{L'}\gamma \cdot f=f'\cdot \partial^-_Ls.$$
\end{defin}
To motivate this, recall the {\em bundle categorification}
of an extension of groups:
$$1\to N\to E\to G\to 1$$

%
If $s,s':G\to E$ are sections,
and $L=C_s, L'=C_{s'}:G \to Aut(N)$ the induced quasi-actions,
then $(s,f)$ and $(s',f')$ are monoidal functors from the
{\em base category} $\C{G}=\C{B}_G$
- using base categorification to avoid constraints for functorial
morphisms due to the presence of morphisms between distinct
elements of $G$ - to the {\em bundle category} $\C{E}=\C{C}(\C{E})$
under bundle categorification. Now $\gamma:s\to s'$,
defined as $\gamma(g)=s'(g)s(g)^{-1}$, is a monoidal morphism:
$$\xymatrix @C=2pc @R=.5pc {
s(ab) \ddto^{f} \rto^{\gamma(ab)} & s'(ab) \ddto^{f'}\\
& & G \ar@/^1pc/[rr]^{(s,f)}  \ar@{}[rr]|{\Downarrow \gamma}
    \ar@/_1pc/[rr]_{(s',f')} & & E \\
s(a)\otimes s(b) \rto^{\gamma_a \otimes \gamma_b} &
  s'(a)\otimes s'(b)
}$$
where $\otimes$ is the monoidal product in $\C{E}$. Recall that
the monoidal product preserves vectors (corollary \ref{cor:tens}),
so that all the above morphisms are vectors, and that all vector diagrams
commute.
%
%
\begin{th}\label{T:monstr}
Let $s,s':E\to G$ two sections of the group extension $1\to N\to E\to G\to 1$
and $\C{F}_N,\C{C}_E,\C{B}_G$ the strict monoidal categories
associated through fiber, bundle and base categorification. Then:

(i) $(s,f)$ and $(s',f')$ are monoidal functors, where
$f=\delta_Ls$ and $f'=\delta_{L'}s'$.

(ii) $\gamma:s\to s'$ defined by $\gamma(g)=s'(g)s(g)^{-1}$
is a monoidal morphism:
$$\diagram
s(ab) \dto_{f_{a,b}} \rto^{\partial^-_L\gamma_{a,b}}
  & s'(ab) \dto^{f'_{a,b}}
 & \partial^-s \dto_{f} \rto^{\partial^-\gamma} & \dto^{f'} \partial^-s'\\
s(a)s(b) \rto_{\partial^+_{L'}\gamma_{a,b}} & s'(a)s'(b)
 & \partial^+s \rto^{\partial^+\gamma} & \partial^+s' 
\enddiagram$$
where $\partial^\pm$ are the coboundary maps of the
categorical parity complex associated to the strict monoidal category
$(\C{F}_N,\otimes)$:
$$\partial^+\gamma_{a,b}=\gamma_a\otimes\gamma_b:
s(a)s(b)\overset{\partial^+_{L'}\gamma(a,b)}{\To} s'(a)s'(b)$$
$$\partial^-\gamma_{a,b}=\gamma(a\otimes b):
s(ab)\overset{\partial^-_L\gamma(a,b)}{\To} s'(ab)$$
The above morphisms are vectors and the underlying
elements of $N$ are given by the
coboundary maps $\partial^-_L$ and $\partial^+_{L'}$,
of the parity quasicomplex $C^\bullet(G,N)$.
\end{th}
\begin{pf}
(i) $f$ and $f'$ are monoidal structures by proposition \ref{P:functorf}.

(ii) We only need to prove the second part of the statement.
Interpreting $\gamma$ as a monoidal morphism, we have:
$$\partial^+\gamma_{a,b}=\gamma_a\otimes \gamma_b, \qquad
  \partial^-\gamma_{a,b}=\gamma_{ab}$$
while as a group cochain:
$$\partial^+_{L'}(a,b)=L'_a(\gamma(b))\gamma(a), \qquad
  \partial^-_L(a,b)=\gamma(ab).$$
The element of $N$ underlying the map $\gamma_a\otimes \gamma_b$
- a vector - is $s'(a)s'(b)(s(a)s(b))^{-1}$ and equals
$\partial^+_{L'}\gamma(a,b)$:
$$ s'(a)\gamma(b)s'(a)^{-1}\gamma(a)=
  s'(a)\ (s'(b)s(b)^{-1})\  s'(a)^{-1}\ (s'(a)s(a)^{-1})=
  s'(a)s'(b)s(b)^{-1}s(a)^{-1}.$$
\end{pf}
As a consequence, since monoidality is preserved under conjugation
(see \cite{Sa,I}), 
a 2-cocycle (monoidal structure) conjugated by a 1-cochain
(monoidal morphism) yields another 2-cocycle, provided they
are represented by group extensions (see section \ref{SS:ext}).
%
\begin{cor}\label{C:Z2}

(i) The set of absolute integrable 2-cocycles $Z^2(G,N)$ is stable
under conjugation by $C^1(G,N)$.

(ii) Any 2-cochain cohomologous to an absolute integrable 2-cocycle
is necessary a cocycle.
\end{cor}

For 1-cochains, we have:
\begin{lem}\label{L:aff1}
Let $(s,L=C_s\circ L^0)$ and $(s',L'=C_{s'}\circ L^0)$ be cochains
corresponding to the same outer quasiaction $[L^0]$, and define
$(\gamma,L')$ by $\gamma=s's^{-1}$:
$$ (s,L)\overset{(\gamma,L')}{\To}(s',L').$$
Then, $(\gamma,L')$ is a cocycle iff
$\bar{\delta}(s,L)=\bar{\delta}(s',L')$.

In particular, if $(s,L)$ is a cocycle, then $(\gamma,L')$ is a cocycle
iff $(s',L')$ is a cocycle.
\end{lem}
\begin{pf}
Note that $\partial^+_Ls(a,b)=s(a)L^0_a(s(b))s(a)^{-1}s(a)=s(a)L^0_a(s(b))$,
and $\partial^+_{L'}s'(a,b)=s'(a)L^0_a(s'(b))$. Since:
\begin{align}
\partial^+_{L'}\gamma(a,b)&=L'_a(\gamma(b))\gamma(a)=
  s'(a)L^0_a(s'(b))L^0_a(s(b))^{-1}s'(a)^{-1}s'(a)s(a)^{-1}\notag\\
  &=s'(a)L^0_a(s'(b))\ (s(a)L^0_a(s(b)))^{-1}\notag\\
  &=(\partial^+_{L'}s')(\partial^+_Ls)^{-1}(a,b)
\end{align}
(see Remark \ref{R:order})
$$\partial^-_{L'}\gamma(a,b)=\gamma(ab)=s'(ab)s(ab)^{-1}
  =(\partial^-_{L'}s')(\partial^-_Ls)^{-1}(a,b)$$
$(\gamma,L')$ is a cocycle iff:
$$(\partial^+_{L'}s')(\partial^+_Ls)^{-1}=
  (\partial^-_{L'}s')(\partial^-_Ls)^{-1}$$
$$(\partial^-_{L'}s')^{-1}(\partial^+_{L'}s')=
  (\partial^-_Ls)^{-1}(\partial^+_Ls)^{-1}$$
i.e. $\bar{\delta}(s,L)=\bar{\delta}(s',L')$. Note that, in general,
$\bar{\delta}_Lc=1$ iff $\delta_Lc=1$.
\end{pf}
\begin{rem}
Note that standard cocycles $(s,C_s)$ ($[L^0]=1$) are the usual morphisms of
groups, since $\partial^+_{C_s}s(a,b)=(s(a)s(b)s(a)^{-1})s(a)=s(a)s(b)$.

Thus to a pair of morphisms of groups, one may associate another
cocycle:
$$\partial((s,C_s),(s',C_{s'}))=(s's^{-1},C_{s'}).$$
An affine group structure will be investigated later.

Also note that, when $s$ and $s'$ are standard cocycles,
the corresponding (trivial) 2-cocycles $(1,L)$ and $(1,L')$
are week cohomologous through 1-cocycles. In particular, since $L=L'$
iff $s$ and $s'$ differ by a central 1-cochain $\gamma\in C^1(G,Cen(N))$,
the isotropy group of $(1,L)$ is $Z^0(G,Cen(N))$.
\end{rem}


\begin{defin}
The parity quasicomplex $(C^\bullet(G,N),\partial^\pm)$ is called
the {\em standard parity quasicomplex} associated to the pair
$(G,N)$ of groups.

The corresponding cohomology spaces - 
classes of cohomologous cocycles (see section \ref{S:pqc}) -
are called {\em total cohomology spaces} and denoted as $H^k_{nc}(G,N)$.

Corresponding to the relation $\underset{wk}{\sim}$ - weak cohomologous -
we define the {\em weak cohomologous spaces} $H^k_{wk}(G,N)$, of $G$
with coefficients in $N$.

Fix a left quasiaction $L$ of $G$ on $N$. The cohomology spaces
$H^k(G, _LN)$, of $G$ with coefficients in $N$, correspond to $C^k(G, _LN)$,
the $L$-component of the total parity quasicomplex.
\end{defin}
Obviously $H^k_{nc}(G,N)=\coprod_L H^k(G, _LN)$.

\subsection{Integrable Cochains and Extensions}\label{SS:intc}

\begin{lem} If $(s,L)$ is an integrable 1-cochain, $\C{L}=C_L$
and $f=\delta_Ls$, then:

(i) $(f,L)$ is an integrable 2-cocycle, where $f=\delta_Ls$.

(ii) $(L,\C{L})$ is an integrable standard 1-cochain, and
$(F,\C{L})$ is an integrable 2-cocycle, where $F=\delta_\C{L}L$.
\end{lem}
\begin{pf}
(i)
Equation (MC) holds on $\C{I}$ since $(s,L)$ is integrable.
Then $\delta^2_L\circ\delta^1_Ls=1$ by lemma \ref{L:dds}.

(ii) By (ii) proposition \ref{P:int}, a standard 1-cochain
is absolute integrable, and equation (MC) holds globally for $(F,\C{L})$.
Moreover, by (ii) lemma \ref{L:dds}, $\delta^2_\C{L}L=1$
and $F$ is a 2-cocycle.
\end{pf}
Any integrable 2-cocycle $(f,L)$ defines an extension
$E_{f,L}=\C{I}\times G$, where $\C{I}\subset N$ is the normal subgroup
generated by the orbit of $s$ under the quasiaction $L$. The group
structure is given by:
\begin{equation}\label{E:mult}
(n_1,g_1)\cdot(n_2,g_2)=(n_1L_{g_1}(n_2)f(g_1,g_2),g_1g_2), \quad n_i\in N,
g_i\in G
\end{equation}
If $(f,L)$ is absolute integrable, i.e. equation (MC) holds on $N$,
then as usual (see \cite{Brown}, section 4.6) $G$ extends by $N$:
$E_{f,L}=N\times G$.
Moreover, if $\tilde{s}(g)=(1,g)$ is the {\em canonical section} of
$E_{f,L}$, then:
\begin{prop}
$1\to N\to E_{f,L}\to G\to 1$ is a group extension, 
$L=C_{\tilde{s}}$ and $f=\delta_L{\tilde{s}}$.
\end{prop}
A more difficult result for 2-cocycles (theorem \ref{T:ext}) is proved in
section \ref{SS:ext}. If the 2-cocycle is absolute integrable, i.e. if
equation (MC) holds globally on $N$, then by corollary \ref{C:assoc}
the group operation \ref{E:mult} is associative and the above proposition
follows.

If we start with an integrable 1-cochain $(s,L)$ and consider the extension
corresponding to the integrable 2-cocycle $f=\delta_Ls$, the section
$\tilde{s}(g)=(s(g),g)$ and the canonical section $s_c$, then, as before
we have $L=C_{s_c}$, $f=\delta_Ls_c$.
Moreover:
\begin{lem}
In the above context, let $L'=C_{\tilde{s}}$ and $f'=\delta_{L'}\tilde{s}$.
Then:

(i) $L'=C_s\circ L$ and $(s,L)\overset{(\gamma,C_\gamma)}{\to}(1,L')$
define the same outeraction $[L]=[L']$,
where $\gamma=s^{-1}$.

(ii) With $L'=C_s$ and $f'=\delta_{L'}s$, $(f',L')$ is cohomologous to
$(f,L)$.

(iii) The extensions $E_{f,L}$ and $E_{f',L'}$ are isomorphic.
\end{lem}
%

We will give a categorical interpretation for the equation (MC).
\begin{th}\label{T:MC}
Let $(f,L)$ be a 2-cochain. If $(f,L)$ is an absolute integrable 2-cocycle,
then:

1) $\C{E}=\C{C}(N\overset{f,L}{\times}G)$ is a strict monoidal category and
$(\tilde{s},f):\C{B}_G\to \C{E}$ is a monoidal
functor, where $\tilde{s}$ is the canonical section.

2) $\C{A}=\C{C}(Aut(N)\overset{L,\C{L}}{\times}G)$ is a strict monoidal
category and $(\tilde{L},F):\C{B}_G\to \C{A}$ is a monoidal functor,
where $\C{L}=C_L$ and $F=\delta_\C{L}L$.

3) $(\B{L},C\times id):N\times G \to Aut(N)\times G$ is a monoidal functor
and the following diagram commutes:
$$\diagram
\C{B}_G \rto^{\tilde{s}} \drto_{\tilde{L}} &
  \C{C}(N\overset{f,L}{\times}G) \dto^{\B{L}}\\
& \C{C}(Aut(N)\overset{L,\C{L}}{\times}G)
\enddiagram$$
where $\B{L}(n,g)=(L_g,g), n\in N, g\in G$.
\end{th}
\begin{pf}
If $f$ is an absolute integrable 2-cocycle, then $E_{f,L}$ is an extension
of $G$ by $N$. By categorification $\C{C}(E_{f,L})$ is a strict monoidal
category, $s$ is a functor and $f$ a monoidal structure.

$(L,\C{L})$ is a standard 1-cochain, thus absolute integrable, and
$(F,\C{L})$ is an absolute integrable 2-cocycle. By categorification,
$\C{A}$ is a strict monoidal category, $\tilde{L}$ a functor and
$F$ a monoidal structure.

Note that the equation (MC) is equivalent to $F=C\circ f$, which
gives the relation between the monoidal structures.
Also $\tilde{L}=(\B{L}, id)\circ \tilde{s}$.

What is left to check is that $(C\times id)$ is globally a monoidal
structure on $\C{E}$.
\end{pf}
The above statement holds for integrable cochains, if $N$ is replaced by
the corresponding subgroup $\C{I}(f,L)$.



%

\begin{rem}
Denote by $Aut^k(N)$ the group obtained by taking iteratively the
automorphisms group: $Aut(Aut(...(N)...)$, with the convention
$Aut^0(N)=N$. Then we have a double quasi-complex
$(C^p(G,Aut^q(N)), \delta^p_q, C_*)$, since in general $C^2\ne 1$.
%
%
%
The maps are strict morphisms - in the
{\em vertical direction} $q$ - and it is a complex in the
{\em horizontal direction}.
If $C^K=1$ (for example if $Aut^K(N)=1$), then it is a $K$-complex
\cite{Ka} in the vertical direction.
\end{rem}

\begin{rem}
If $N$ is abelian, $C_*=0$ and
the complexes $C^\bullet(G,Aut^q(N))$ are not coupled.

Note also that the orbit $C_s\circ L$  of a quasi-action $L$ is one
point.
\end{rem}


\subsection{Quasi-extensions and $H^3$}\label{SS:qext}

Let $(f,L)$ an arbitrary 2-cochain. Define a multiplication on $E=N\times G$,
not necessary associative, by equation:
\begin{equation}\label{E:mult1}
(n_1,g_1)\cdot(n_2,g_2)=(n_1L_{g_1}(n_2)f(g_1,g_2),g_1g_2), \quad n_i\in N,
g_i\in G
\end{equation}
Denote by $(n,g)^*=(m,h)$ the right inverse:
\begin{equation}
(n,g)\cdot(m,h)=(1,1) \quad \text{iff}\quad
  h=g^{-1},\ L_g(m)f(g,g^{-1})=n^{-1}
\end{equation}
Consider the category $\C{E}$, as in bundle categorification. Objects
are elements $(n,g)$ and the only morphisms are $Hom((n,g),(n',g))=N$,
with composition the multiplication in $N$.
Define a product $\otimes$ in $\C{E}$ as before:
\begin{gather}
\diagram
(n_1,g) \dto^{k \quad\otimes} &
  e \dto^{n \qquad =\quad} & (n_1,g)e \dto^{(n_2,g) (n,1) (n_1,g)^*}\\
(n_2,g) &  e' & (n_2,g) e' &
\enddiagram\label{E:red1}\end{gather}
A direct computation shows that it is well defined:
$$((n_2,g)(g,1))(n_1,g)^*=(n_1,g)((g,1)(n_1,g)^*)=(n_2L_g(n)n_1^{-1},1).$$
$(\C{E},\otimes)$ is a category with product. Any central 3-cochain
$\alpha\in C^3(E,Cen(N))$ defines a quasi-associator (non-coherent).

Define the {\em vector} from
$(n_1,g)$ to $(n_2,g)$ as the morphism corresponding to the $N$ element:
$$(n_2,g)\cdot(n_1,g)^*=(n_2L_g(m)f(g,g^{-1}),1)=(n_2n_1^{-1},1).$$
It follows that composition preserves vectors, since a direct check
shows $(e_1e_2^*)(e_2e_3^*)=e_1e_3^*$, where $e_i=(n_i,g)$.
Denote by $\C{E}_r$ the
subcategory of $\C{E}$ obtained by restricting morphisms to vectors.
Each connected component is simply connected:
any diagram in $\C{E}_r$ commutes.

It is easy to see that the product $\otimes$ respects vectors. If in
equation \ref{E:red1} $e=(n_1',g')$ and $e'=(n_2',g')$ then
$n=n_2'(n_1')^{-1}$ and $(n_1,g)e=(n_1L_g(n_1')f(g,g')=x,gg')$,
$(n_2,g)e'=(n_2L_g(n_2')f(g,g')=y,gg')$. Then
$(y,gg')(x,gg')^*=(yx^{-1},1)$ and:
$$yx^{-1}=n_2L_g(n_2')L_g(n_1')^{-1}n_1^{-1}=n_2L_g(n)n_1^{-1}.$$
Thus
$(\C{E}_r,\otimes)$ is a subcategory with the induced product. Define
$\tilde{\alpha}\in C^3(E,N)$ as being given by the unique vector
$(e_1e_2)e_3\to e_1(e_2e_3)$, for any $e_i\in E$.
Then $\tilde{\alpha}$ is natural - any diagram in $\C{E}_r$ commutes -
and $\delta{\tilde{\alpha}}=1$.
Thus $(\C{E}_r,\otimes,\tilde{\alpha})$ is a monoidal category.

Identifying $N$ with its image through the canonical embedding
$j(n)=(n,1)$ (and $Aut(N)$ with $Aut(j(N))$), we have:
\begin{lem}\label{L:qext}
(i) $N$ is normal in $(E,\cdot)$, i.e.
$$C_e(n)=(e (n,1)) e^*=e((n,1)e^*)\in N$$
and $C_{(n,g)}=C_n\circ L$. In particular $L=C_s$, where
$s$ is the canonical section $s(g)=(1,g)$.

(ii) With $\tilde{L}=C_s:E\to Aut(N)$, $(\tilde{\alpha},\tilde{L})$
is a 3-cocycle in $C^3(E,N)$.

(iii) The vector $s(gg') \to s(g)s(g')$ is $f(g,g')$.

(iv) The natural projection $\pi:E\to G$ is compatible with the products,
i.e. $\pi(e_1e_2)=\pi(e_1)\pi(e_2)$.
\end{lem}
Note first that the conjugation of $E$ is well defined only on $N$, so
that $\partial^+_Ls$ is undefined.
\begin{pf}
(ii) was proved above.

A direct computation shows (i),(iii) and (iv) hold.
\end{pf}
Consider the strict monoidal categories $\C{B}_G$ and $\C{F}_N$,
of $G$ and $N$ corresponding to base and fiber categorification.
%
%
\begin{th}\label{T:qext}
$(\C{E}_r,\otimes,\tilde{\alpha})$ is a monoidal category.
$j$ and $\pi$ are strict monoidal functors.
$(s,f)$ is a monoidal functor.
\end{th}
\begin{defin}\label{defin:qext}
$1\to N\to E\to G\to 1$ is called the {\em quasi-extension} associated
to the (normalized) 2-cochain $(f,L)$, and {\em represented}
by $\tilde{s}$, the canonical section:
$$\diagram
E \dto_{C} & G \lto_{\tilde{s}} \dlto^{L}\\
Aut(N)
\enddiagram$$
$(\C{E}_r(f,L),\otimes,\tilde{\alpha})$ is the associated monoidal
category through {\em reduced bundle categorification}.
\end{defin}
In order for $\tilde{\alpha}$ to depend only on the base elements
$g\in G$, we need additional assumptions.
\begin{prop}
If $(f,L)$ is an integrable 2-cochain and $\alpha=\delta_Lf$,
then $\tilde{\alpha}$ factors
through $\pi$ iff $\alpha=\delta_Lf$ is valued in the centralizer
of $\C{I}(f,L)$.
$$\diagram
E^3 \dto_{\pi} \rto^{\tilde{\alpha}} & N\\
G^3\urto_{\alpha}
\enddiagram$$
In this case, $\delta_L\alpha=1$.
\end{prop}
\begin{pf}
By (i) lemma \ref{L:diagr}, 
$\tilde{\alpha}_{e_1,e_2,e_3}=\alpha_{g_1,g_2,g_3}$
iff
\begin{align}\label{E:alfa}
\alpha_{g_1,g_2,g_3}
n_1L_{g_1}(n_2)f(g_1,g_2)&L_{g_1g_2}(n_3)f(g_1g_2,g_3)\\
&=n_1L_{g_1}(n_2)L_{g_1}(L_{g_2}(n_3))L_{g_1}(f(g_2,g_3))f(g_1,g_2g_3)
  \notag
\end{align}
where $n_i\in \C{I}$. 
Recall that $\alpha=\delta_Lf$ verifies:
\begin{equation}\label{E:co}
L_{g_1}(f(g_2,g_3))f(g_1,g_2g_3)=
  \alpha_{g_1,g_2,g_3} f(g_1,g_2) f(g_1g_2,g_3).
\end{equation}
If (MC) equation holds, then equation \ref{E:alfa} becomes:
\begin{align}
\alpha_{g_1,g_2,g_3}
n_1L_{g_1}(n_2)&L_{g_1}(L_{g_2}(n_3))f(g_1,g_2)f(g_1g_2,g_3)\\
&=n_1L_{g_1}(n_2)L_{g_1}(L_{g_2}(n_3))L_{g_1}(f(g_2,g_3))f(g_1,g_2g_3)
  \notag
\end{align}
or equivalently, by equation \ref{E:co}:
\begin{align}
\alpha_{g_1,g_2,g_3}
n_1L_{g_1}(n_2)&L_{g_1}(L_{g_2}(n_3))f(g_1,g_2)f(g_1g_2,g_3)\\
&=n_1L_{g_1}(n_2)L_{g_1}(L_{g_2}(n_3))
  \alpha_{g_1,g_2,g_3} f(g_1,g_2) f(g_1g_2,g_3).
\end{align}
i.e.:
\begin{align}
\alpha_{g_1,g_2,g_3}
n_1L_{g_1}(n_2)&L_{g_1}(L_{g_2}(n_3))\\
&=n_1L_{g_1}(n_2)L_{g_1}(L_{g_2}(n_3))
  \alpha_{g_1,g_2,g_3}.
\end{align}
Now take $n_2=n_3=1$ to obtain that $\alpha$ commutes with $\C{I}$.

Conversely, if $\alpha$ commutes  with $n_3, f's$ and $L_a(f(b,c)$,
from equation (MC) and \ref{E:co} follows that
$\tilde{\alpha}_{e_1,e_2,e_3}=\alpha_{g_1,g_2,g_3}$, i.e.
$\tilde{\alpha}$ factors through $\pi$.
\end{pf}

\subsection{Group extensions and $H^2$}\label{SS:ext}

Next we will prove that 2-cocycles are integrable. Let $(f,L)$ be
a 2-cocycle, $\alpha=\delta_Lf=1$.
Consider the categorical interpretation. The reduced bundle
category $(\C{E}_r,\otimes,\tilde{\alpha})$
is a monoidal category (theorem \ref{T:qext}),
and $\tilde{\alpha}$ verifies the pentagon diagram:
\begin{gather}
\xymatrix @C=.5pc @R=1pc {
& & (s_as_b)(s_cs_d) \drrto^{\tilde{\alpha}_3}\\
((s_as_b)s_c)s_d \drto^{\tilde{\alpha}_4} \urrto^{\tilde{\alpha}_1} & & & &
  s_a(s_b(s_cs_d)) \\
& (s_a(s_bs_c))s_d \rrto^{\tilde{\alpha}_2} & &
  s_a((s_bs_c)s_d) \urto^{\tilde{\alpha}_0}
}\label{D:penta}\end{gather}
Recall that $G$ is a strict monoidal category and $f$ is a monoidal
structure for the functor $s:G\to E$ - the canonical section.

We will prove that any 2-cocycle is integrable, i.e. the cocycle
condition:
\begin{gather}
\diagram
s_{(ab)c} \ar@{=}[d] \rto^{f(ab,c)} & s_{ab}s_c \rto^{f(a,b)} &
  (s_as_b)s_c \ar@{=}[d]\\
s_{a(bc)} \rto^{f(a,bc)} & s_as_{bc} \rto^{L_a(f(b,c))} &
  s_a(s_bs_c)
\enddiagram\label{D:cocycle}\end{gather}
implies that the equation (MC) holds on the corresponding subgroup $\C{I}$
(compare \cite{Brown}, page 105).

Note first that:
\begin{lem}\label{L:diagr1}
(i) $\tilde{\alpha}_{s(a),s(b),s(c)}=\alpha(a,b,c)=1$ for any $a,b,c\in G$.

(ii) In the diagram \ref{D:penta}, $\alpha_4=1$ and $\alpha_0=1$.
\end{lem}
\begin{pf}
(s,f) is a monoidal functor, $G$ is a strict monoidal category, so that:
$$\tilde{\alpha}_{s_a,s_b,s_c}=\partial^+f \cdot s(1) \cdot \partial^-f=1$$
where $\partial^\pm f$ are the categorical coboundary maps.

Alternativly, by a direct computation:
$$((1,a)(1,b))(1,c)=(f(a,b),ab)(1,c)=(f(a,b)f(ab,c),(ab)c)$$
and
$$(1,a)((1,b)(1,c))=(1,a)(f(b,c),bc)=(L_a(f(b,c))f(a,bc),a(bc)).$$
As a consequence, we have (ii).
\end{pf}
Recall that multiplication of a map $e_1\overset{k}{\to}e_2$ by $I_{(n,g)}$
on the right is trivial, and to the left gives the map
$y:(n,g)e_1 \to (n,g)e_2$ with products for source and target,
and corresponding $y\in N$ (or $(y,1)$ if viewed as an element of $E_r$):
$$(y,1)=(n,g)(k,1)(n,g)^*=(nL_g(k)n^{-1},1).$$
%
\begin{lem}\label{L:diagr}
If $(f,L)$ is a 2-cocycle, then:

(i) If $e_i=(n_i,g_i)$ then $\tilde{\alpha}_{e_1,e_2,e_3}:(x,g_1g_2g_3)
\to (y,g_1g_2g_3)$, where:
$$x=n_1L_{g_1}(n_2)f(g_1,g_2)L_{g_1g_2}(n_3)f(g_1g_2,g_3)$$
$$y=n_1L_{g_1}(n_2)L_{g_1}(L_{g_2}(n_3))L_{g_1}(f(g_2,g_3))f(g_1,g_2g_3)$$

(ii) $\tilde{\alpha}_{s(a)s(b),s(c),s(d)}=1$, so that $\alpha_1=1$.

(iii) $\tilde{\alpha}_{s(a),s(b)s(c),s(d)}=1$, so that $\alpha_2=1$.

(iv) $\tilde{\alpha}_{s(a),s(b),s(c)s(d)}=L_a(L_b(f(c,d))f(a,b)
L_{ab}(f(c,d))^{-1}f(a,b)^{-1}$.
\end{lem}
\begin{pf}
(i) (computation)

(ii) Use (i) with $n_1=f(a,b), n_2=1, n_3=1$ and the cocycle condition.

(iii) Similar, with $n_1=1, n_2=f(b,c), n_3=1$.
$$((1,a)(f(b,c),bc))(1,d)=(L_a(f(b,c))f(a,bc)f(abc,d),abcd)$$
and
$$(1,a)((f(b,c),bc)(1,d))=(L_a(f(b,c))L_a(f(bc,d))f(a,bcd),abcd)$$
are equal since $\alpha_{a,bc,d}=1$.

(iv) Now $n_1=1, n_2=1, n_3=f(c,d)$.
\end{pf}
\begin{cor}\label{C:assoc}
If $(f,L)$ is an absolute integrable 2-cocycle, then the associator
of the monoidal category $\C{E}_r$ is trivial.
\end{cor}
\begin{pf}
Equation (MC) holds for any $n\in N$:
\begin{equation}\label{E:MCexpl}
L_a(L_b(n))f(a,b)=f(a,b)L_{ab}(n), \quad any\ a,b,c\in G.
\end{equation}
Then, by (i) lemma \ref{L:diagr}, $\tilde{\alpha}=1$ iff
$f$ is a 2-cocycle.
\end{pf}
\begin{cor}
All the vertices of the pentagonal diagram \ref{D:penta} are equal,
therefor $\tilde{\alpha}_3=1$.
\end{cor}
\begin{pf}
Since $\tilde{\alpha}$ is the vector between its source and target,
(ii) lemma \ref{L:diagr1} and (ii), (iii) lemma \ref{L:diagr}
imply that all vertices are the same.
\end{pf}

Although $\C{I}=<f,L>$, it is not necesary that $\{s(g)=(1,g)|g\in G\}$
generates entire $\C{E}_r=\C{I}\times G$.

In any case $\tilde{\alpha}\cong 1$, and $\C{E}_r$ is a strict monoidal
category, as a consequence of the following proposition.
\begin{prop}\label{p:mc}
The following are equivalent:

(i) $\tilde{\alpha}_3=1$.

(ii) 
\begin{equation}\label{E:mc}
L_a(L_b(f(c,d)))f(a,b)=f(a,b)L_{ab}(f(c,d)), \quad any\ a,b,c\in G
\end{equation}

(iii) Equation (MC) holds on $\C{I}(f,L)$.
\end{prop}
\begin{pf}
(i) and (ii) are equivalent, by (iv) lemma \ref{L:diagr}.

Obviously (MC) implies (ii).

Note first that equation \ref{E:mc} for the element $x=f(c,d)$ implies
equation \ref{E:mc} for $x^{-1}$. By taking inverses:
$$f(a,b)^{-1}L_a(L_b(f(c,d)^{-1}))=L_{ab}(f(c,d)^{-1})f(a,b)^{-1}$$
and then multiply both sides, appropriately by $f(a,b)$.

Elements generated by applying $L_a$ can be expressed as products
of elements of the form $x=f(c,d)$ by using cocycle condition
\ref{D:cocycle} for $f$. For example:
\begin{align}
L_a(L_b(L_h(c,d)))f(a,b)&=L_a(L_b(x_2x_1y^{-1}))f(a,b)\notag\\
 &=L_a(L_b(x_2))L_a(L_b(x_1))L_a(L_b(y^{-1}))f(a,b)\notag
\end{align}
Now $f(a,b)$ ``passes through'' $L_aL_b$ yielding $L_{ab}$. Combining
again the factors - $L_g$ are morphisms of groups - gives equation \ref{E:mc}
for $L_h(c,d)$.
\end{pf}
As a consequence, using (i) from lemma \ref{L:diagr},
$\tilde{\alpha}=1$ is trivial, and we have:
%
%
\begin{th}\label{T:ext}
Any p-cocycle $(c,L)\in C^p(G,N)$ is integrable, where $0\le p\le 2$. 

If $(f,L)$ is a 2-cocycle, then
$E_{f,L}=\C{I}(f,L)\overset{(f,L)}{\rtimes}G$ is a group extension of
$G$ by the subgroup $\C{I}(f,L)$ of $N$ and $(f,L)$ is represented
by the 1-cochain $(s,L)\in C^1(G,E)$, where $s$ is the canonical section:
$$ L=C_s \qquad f=\delta_Ls.$$

The associated monoidal category $(\C{E}_r,\otimes,\tilde{\alpha})$ is
strict and $(s,\delta s:\partial^-s\to \partial^+s):G\to E$ is a
monoidal functor.
\end{th}
Thus 2-cocycles are stratified by the subgroups of $N$.

The absolute integrable 2-cocycles, for which equation (MC) holds
globally have a nice moduli space ($H^2(G,Cen(N)$; see \cite{Brown}).
The vector $f_2f_1^{-1}$ for two such cocycles is central.
Recall (\cite{Brown}, section 4.6)
that the set $\C{E}(G,N,\psi)$ of isomorphism classes of extensions
inducing a given outeraction $\psi:G\to Out(N)$ is in bijection with
$H^2(G,Cen(N))$, provided that a certain obstruction $u\in H^3(G,Cen(N))$
(associator) corresponding to the exact
sequence \ref{D:cen} vanishes, or else it is empty.

%
Still, a 2-cocycle need not be absolute integrable.
\begin{defin}
A cocycle $(f,L)$ is called {\em irreducible} iff $\C{I}(f,L)=N$.
\end{defin}
\begin{cor}
(i) An irreducible 2-cocycle is absolute integrable.

(ii) To any irreducible 2-cocycle corresponds an extensions of $G$ by $N$.
\end{cor}

The relation between irreducible
2-cocycles, central 2-cocycles ($H^2(G,Cen(N))$) and stability
will be studied later.

Conversely, let
\begin{equation}
1\to N\overset{i}{\to} E\overset{\pi}{\to} G\to 1\tag{$\C{E}$}
\end{equation}
be an extension, where we may assume $i$ being the inclusion of $N$ in $E$,
to simplify the notation. 
Any section $s:G\to E$ of $\pi:E\to G$ induces a quasiaction
$L=C_s:G\to Aut(N)$, of $G$ on $N$, with a unique {\em outer quasiaction}
$\Psi=[L]:G\to Out(N)$, where $Out(N)=Aut(N)/Int(N)$ is the quotient
modulo inner actions.

Recall that $f=\delta_Ls:G\times G\to N$ is a 2-cocycle relative to $L$.
Moreover $(f,L=C_s)$ is absolute integrable, but not necessary irreducible. 

%
\subsection{Split extensions and $H^1$}\label{SS:split}

For split extensions with abelian kernel see \cite{Brown}.

Assume now that the 2-cocycle $(f,L)$ is
cohomologous to $(1,L)$
$(\partial^+_L\gamma)\cdot 1=f\cdot (\partial^-_L\gamma)$,
i.e. $f$ is the coboundary of $\gamma$. Represent $(f,L)$ as in
section \ref{SS:ext} by $(s,L)$, where $s$ is the canonical section of the
associated extension $E=\C{I}\overset{(f,L)}{\times}G$, where $\C{I}$
is the holonomy group of $(s,L)$.

$\C{E}$ {\em splits} iff there is a section $s'$ of $\pi$ which is a morphism
of groups, i.e. if $s'$ itself is a 1-cocycle for the corresponding
quasiaction $L'=C_{s'}$.

Define $s'(g)=(\gamma(g)^{-1},g)$ - so that $\gamma:s'\to s$ -,
and $L'=C_{s'}$.
Consider $E'=\C{I}\overset{(1,L')}{\times}G$ with multiplication:
$$(n_1,g_1)(n_2,g_2)=(n_1L'_{g_1}(n_2),g_1g_2)$$
and $\phi:E'\to E$ defined by $\phi(n,g)=(n\gamma(g)^{-1} ,g)$.
\begin{prop}\label{P:split}
(i) $\delta_{L'}s'=1$ and $s':G\to E$ is a morphism of groups.

(ii) $\phi(e_1e_2)=\phi(e_1)\phi(e_2)$.

Thus $L'$ is a morphism and $E'$ is a semi-direct product
isomorphic to $E$.
$$\diagram
\C{I} \overset{L'}{\rtimes} G \rto^{\phi} &
  E=\C{I}\overset{(f,L)}{\times}G \ar@/_/[d]_{\pi}\\
& G \ulto^{s}
\ar@/_/[u]_{s'=(\gamma^{-1},id)} 
\enddiagram$$
\end{prop}
\begin{pf}
Since $L'=C_{s'}$, with multiplications in $E$, we have:
\begin{align}
\delta_{L'}s'(a,b)&=s'(a)s'(b)s'(ab)^{-1}\notag\\
  &=(\gamma(a)^{-1},a)(\gamma(b)^{-1},b)(\gamma(ab)^{-1},ab)^{-1}\notag\\
  &=(\gamma(a)^{-1}L_a(\gamma(b)^{-1})f(a,b),ab)(m,(ab)^{-1})\notag\\
  &=(\gamma(a)^{-1}L_a(\gamma(b)^{-1})f(a,b)L_{ab}(m)f(ab,(ab)^{-1}),1)
  \notag\\
  &=(\gamma(a)^{-1}L_a(\gamma(b)^{-1})L_a(\gamma(b))\gamma(a)
      \gamma(ab)^{-1}\gamma(ab),1)=(1,1)\notag
\end{align}
and $s'$ is a morphism of groups.

A similar computation:
\begin{align}
\phi((n_1,g_1)(n_2,g_2))&=\phi(n_1L'_{g_1}(n_2),g_1g_2)=\notag\\
&=(n_1\gamma(g_1)^{-1}L_{g_1}(n_2)\gamma(g_1)\gamma(g_1g_2)^{-1},g_1g_2)
\notag\\
\phi(n_1,g_1)\phi(n_2,g_2)&=
  (n_1\gamma(g_1)^{-1}L_{g_1}(n_2\gamma(g_2)^{-1})f(g_1,g_2),g_1g_2)
  \notag\\
  &=\phi((n_1,g_1)(n_2,g_2))\notag
\end{align}
shows that the multiplication on $E'$ is carried bijectively on the
group structure of $E$.
\end{pf}
\begin{rem}
In a fixed ``stratum'' ($\C{I}=H$ subgroup of $N$),
2-coboundaries $(f,L)\sim(1,L')$ give split extensions $E'$
isomorphic to trivial 2-coboundaries $(1,L')$, with $(1,1)$ corresponding
to the direct product $N\times G$.
\end{rem}
%
Conversely, let
\begin{equation}
1\to N\overset{i}{\to} E\overset{\pi}{\to} G\to 1\tag{$\C{E}$}
\end{equation}
be a split extension, where we may assume $i$ being the inclusion of $N$
in $E$, to simplify the notation. 

Fix a splitting $s$.
Now $L=C_s$ is an action, and
the pair $(L, \C{L})$ is a 1-cocycle in $C^\bullet(G, Aut(N))$,
where $\C{L}=C_{L}$, so that the associated 2-cocycle
$(f=\delta_Ls,L)=(1,L)$ is absolute integrable.

Note that any two such actions induce the same outer-action $[L]$.

Define the semi-direct product $E'=N \overset{L}{\rtimes} G$.
Then, as usual, $E$ is isomorphic to $E'$
through the bijection $\phi:N\times G\to E$
associated to the splitting morphism $s$, and
defined by $\phi(n,g)=ns(g)$:
\begin{align}
\phi((n_1,g_1)\cdot(n_2,g_2))&=\phi(n_1L_{g_1}(n_2), g_1g_2)\notag\\
&=n_1L_{g_1}(n_2)s(g_1g_2)=n_1s(g_1)n_2s(g_1)^{-1}s(g_1g_2)\notag\\
&=n_1s(g_1)n_2s(g_2)=\phi(n_1,g_1)\phi(n_2,g_2)\notag
\end{align}
The set of sections $s':G\to E$ bijectively correspond with
1@-cochains $\gamma:G\to N$, defined by
$\tilde{\gamma}:(s',L')\to (s,L)$, i.e.
$\gamma:G\to N, \tilde{\gamma}(g)=s(g)s'(g)^{-1}=(\gamma(g),1)$.
Then, by lemma \ref{L:aff1} applied to $C^1(G,E)$,
$(\tilde{\gamma},L)$ is a 1-cocycle iff $(s',L')$
is a 1-cocycle, i.e. if $s'$ is a morphism of groups.
Now $\tilde{\gamma}$ is
$N$-valued, so $s'$ is a morphism of groups iff $(\gamma,L)$ is a cocycle.
$$\diagram
N \overset{L}{\rtimes} G \rto^{\phi} & E \\
& G \ulto^{s'=(\gamma^{-1},id)} \uto_{s'} 
\enddiagram$$
1-cocycles $(\gamma,L)$:
$$L_a(\gamma(b))\gamma(a)=\gamma(ab)
\qquad (\text{additive\ notation}\ \gamma(ab)=\gamma(a)+L_a(\gamma(b))
$$
are also called crossed homomorphisms (see \cite{Brown}, p.88).

Let $(\gamma_i,L)$ be two cohomologous 1-cocycles corresponding to the
morphisms $s_i:G\to E$, i.e. $\partial^+_Ln\gamma_1=\gamma_2\partial^-_Ln$,
with $n\in N$ and $L_g(n)\gamma_1(g)=\gamma_2(g)n$. Equivalently:
$$(\gamma_1(g)^{-1},g)(n^{-1},1)=(n^{-1},1)(\gamma_2(g)^{-1},g)$$
by taking inverses and considering the equality in
$N\overset{L}{\times} G$. Applying $\phi$ we get $s_1(g)n^{-1}=n^{-1}s_2(g)$,
i.e. $s_2=C_n(s_1)$ are $N$@-conjugate.
Thus, as in the abelian case (\cite{Brown}), we have:
%
%
\begin{th}\label{T:split}
For any $G$-group $N$, the $N$-conjugacy classes of splittings
of the split extension:
$$1\to N\to N\overset{L}{\rtimes}G \to G\to 1$$
are in 1:1 correspondence with the elements of $H^1(G, _LN)$.
\end{th}

Note that if $s_1$ and $s_2$ induce the same action $L_1=L_2$, they differ
by a central cocycle $\gamma\in C^1(G,Cen(N))$.

\begin{rem}
The semi-direct products $N \overset{L'}\rtimes G$ belonging to the
connected groupoid $[\C{E}]$ - the isomorphism class of the
split extension $\C{E}$ - are in bijection with $Z^{1}(G, _LN)$ modulo
central cocycles $Z^{1}(G,Cen(N))$, a base point being given
by a splitting $s$. The morphisms $E \to
N \overset{L}\rtimes G$ are parametrized by $Z^{1}(G,Cen(N))$.
\end{rem}




\end{document}